\documentclass[review,onefignum,onetabnum]
{siamart171218}

\usepackage{amsopn}

\usepackage{mathtools}

\usepackage{lipsum}
\usepackage{amsfonts}
\usepackage{graphicx}
\usepackage{epstopdf}
\usepackage{algorithmic}
\usepackage{mathtools}
\ifpdf
  \DeclareGraphicsExtensions{.eps,.pdf,.png,.jpg}
\else
  \DeclareGraphicsExtensions{.eps}
\fi

\newsiamremark{remark}{Remark}
\newsiamremark{hypothesis}{Hypothesis}
\crefname{hypothesis}{Hypothesis}{Hypotheses}
\newsiamthm{claim}{Claim}

\headers{Bilevel learning via recycling Krylov subspaces}{M. J. Ehrhardt, S. Gazzola, and S. J. Scott}

\title{Efficient gradient-based methods for bilevel learning via recycling Krylov subspaces\thanks{
\funding{
MJE acknowledges support from the EPSRC (EP/S026045/1, EP/T026693/1, EP/V026259/1, EP/Y037286/1) and the European Union Horizon 2020 research and innovation programme under the Marie Skodowska-Curie grant agreement REMODEL.
SJS was supported by a scholarship from the EPSRC Centre for Doctoral Training in Statistical Applied Mathematics at Bath (SAMBa), under the project EP/S022945/1 and acknowledges funding by the German Ministry of Science and Technology (BMBF) under grant agreement No. 01IS24072A (COMFORT).}}}
\author{Matthias J. Ehrhardt\thanks{Department of Mathematical Sciences, University of Bath, UK
.}
\and Silvia Gazzola\footnotemark\thanks{ Department of Mathematics, University of Pisa, Italy.}
\and Sebastian J. Scott\footnotemark[3]\thanks{Institute of Mathematics, University of Würzburg, Germany   (\email{sebastian.scott@uni-wuerzburg.de}). }}

\makeatletter
\newcommand*{\addFileDependency}[1]{
  \typeout{(#1)}
  \@addtofilelist{#1}
  \IfFileExists{#1}{}{\typeout{No file #1.}}
}
\makeatother

\newcommand*{\myexternaldocument}[1]{%
    \externaldocument{#1}%
    \addFileDependency{#1.tex}%
    \addFileDependency{#1.aux}%
}

\ifpdf
\hypersetup{
	pdftitle={Efficient gradient-based methods for bilevel learning via recycling Krylov subspaces},
	pdfauthor={M. J. Ehrhardt, S. Gazzola, S. J. Scott}
}
\fi

\myexternaldocument{ex_supplement}

\newcommand\ulfun{\mathcal L} 
\newcommand\ulfunx{\ell} 
\newcommand\llfun{\Phi} 

\DeclareMathOperator{\diag}{diag}
\newcommand\recon{\hat x }
\newcommand\reconp{\hat x(\theta) }
\newcommand\reconpk{\hat x(\theta^{(i)}) }

\newcommand\xtrue{x^\star} 
\newcommand\ynoise{y} 
\newcommand\op{A} 
\newcommand\kry{\mathcal K} 
\newcommand{\grad}{\nabla}

\newcommand{\deriv}{D}
\newcommand{\derivtwo}{D^2}

\newcommand{\jacadj}{ J}

\newcommand\real{\mathbb{R}} 

\DeclarePairedDelimiter{\Vpair}{\|}{\|} 
\DeclarePairedDelimiter{\vpair}{|}{|} 
\newcommand\norm[2][0pt]{\Vpair*{#2}} 
\newcommand\abs[2][0pt]{\vpair*{#2}} 

\DeclareMathOperator*{\argmin}{arg\,min}

\newcommand{\imgsize}{0.2367}

\begin{document}\nolinenumbers
	
	\maketitle
	
	\begin{abstract}
		Many optimization problems require hyperparameters, i.e., parameters that must be pre-specified in advance of their solution, such as regularization parameters and parametric regularizers in variational regularization methods for inverse problems, and dictionaries in compressed sensing. A data-driven approach to determine appropriate hyperparameter values is via a nested optimization framework known as bilevel learning.
		Even when it is possible to employ a gradient-based solver to the bilevel optimization problem, construction of the gradients, known as hypergradients, is computationally challenging,  each one requiring both a solution of a minimization problem and a linear system solve. These systems do not change much during the iterations, which motivates us to apply recycling Krylov subspace methods, wherein information from one linear system solve is re-used to solve the next linear system.
		Existing recycling strategies often employ eigenvector approximations called Ritz vectors. In this work we propose a novel recycling strategy based on a new concept, Ritz generalized singular vectors, which acknowledge the bilevel setting. 
		Additionally, while existing iterative methods primarily terminate according to the residual norm, the newly introduced Ritz generalized singular vectors allow us to define a new stopping criterion that directly approximates the error of the associated hypergradient. The proposed approach is validated through extensive numerical testing in the context of inverse problems in imaging.

	\end{abstract}
	
	\begin{keywords}
		Bilevel learning, recycling Krylov methods, variational regularization, inverse problems
	\end{keywords}
	
	\begin{AMS}
		65F10 
		
	\end{AMS}
	
	\section{Introduction}	\label{sec:intro}
	Large-scale optimization problems that arise in applications such as machine learning  \cite{franceschi2018bph, liu2018ddaa,wang2021fas, wang2022zzo} and inverse problems~\cite{benning2018mrm,engl1996ripa,hansen2010dip,scherzer2009vmi} often require the pre-specification of certain parameters, such as  
	dictionaries for compressed sensing ~\cite{ davidl.donoho2006cs, li2024cdl, mairal2012tdl} and 
	sample weights for a biased dataset~\cite{ ren2018lre,shu2019mle}.
	Referred to as hyperparameters, it is crucial that  appropriate values of these parameters are determined in order to   ensure both convergence of the optimization process and that a meaningful solution is obtained.
	For example, for inverse problems, which arise in many important science and engineering applications such as biomedical, astronomical, and seismic imaging  \cite{arridge2019sip,benning2018mrm, engl1996ripa,   hansen2006dim}, one is tasked to recover some quantity of interest $\xtrue\in\real^n$  given $\op\in\real^{m\times n}$ and indirect measurements $\ynoise\approx \op \xtrue \in\real^m$ thereof.    
	A common approach to determine a meaningful approximation to $\xtrue$ is  to solve an optimization problem such as
	\begin{equation}
		\min_{x\in\real^n}  \left\{ \llfun(x,\theta) =  \norm{\op x - \ynoise}_2^2 + \mathcal R(x,\theta)\right\},\label{eq:var-reg}
	\end{equation}
	where the regularizer $\mathcal R$ is used to encourage a-priori information of the solution in the reconstruction, such as sparsity in some signal domain; see~\cite{benning2018mrm,engl1996ripa,hansen2010dip,scherzer2009vmi}. 
	Here the regularizer is parametrized by some vectors $\theta\in\real^p$ which may represent, for example: the weights for a sum of different regularizers  \cite{delosreyes2016sop}; convolution filters \cite{chen2014rlt,roth2009fe}; parameters of a neural network \cite{amos2017icn,  habring2024nrm, haltmeier2021rip, lunz2018ari}. 
	Determining appropriate values of $\theta$ is crucial to attain a meaningful solution, as a poor choice can lead to a reconstruction that is  dominated by noise  or over-smoothed~\cite{hansen2010dip}.
	
	If one had access to ground-truth and noisy measurements pairs $(\xtrue_k, \ynoise_k),$ which may also be constructed synthetically, one way to determine a data-driven regulariser would be to solve the nested optimization problem
	\begin{equation}
		\min_{\theta\in\real^p} \left\{\ulfun(\theta) = \sum_{k} \frac{1}{2} \norm{\recon_k(\theta) - \xtrue_k}_2^2\right\},\; \text{where }\recon_k(\theta)\text{ solves \eqref{eq:var-reg} with $y=y_k.$}\label{eq:ip-bl}
	\end{equation}
	Then, for a new measurement $y$ from an application similar to those in the training set,  a meaningful solution  can be determined by solving \eqref{eq:var-reg} utilizing said data-driven regularizer.
	
	Due to the learning of hyperparameters and inherent nested structure, \eqref{eq:ip-bl} is called a bilevel learning problem \cite{arridge2019sip,crockett2022bmi,delosreyes2016sop,kunisch2013boa}, sitting in the wider class of bilevel optimization \cite{colson2007obo,sinha2018rbo}.
	The minimization over $x$ in \eqref{eq:var-reg} and $\theta$ in \eqref{eq:ip-bl} are called the lower and upper level problems, respectively.
	Other cost functions in both the  upper and lower problems are possible \cite{deledalle2014sug,fehrenbach2015bid, hintermuller2022dad,santambrogio2024wbl,sixou2021arpa,zhang2020bns}, such that we need not be in a supervised setting as in \eqref{eq:ip-bl}. Indeed, all the following considerations apply to a more general class of bilevel problems than \eqref{eq:ip-bl}, namely, ones which can be solved via a gradient-based method.
	In the case of \eqref{eq:ip-bl}, taking a single training sample for simplicity and assuming enough regularity of the lower level objective function, the so-called hypergradient is given by
	\begin{equation*}
		\nabla\ulfun(\theta) = D^2_{\theta x} \llfun(\recon(\theta),\theta) ^T w_\star,
	\end{equation*}
	where $D_{\theta x}$ denotes the second order mixed derivative and $w_\star$ solves
	\begin{equation}
		D^2_{xx} \llfun(\recon(\theta),\theta) \;w_\star = \recon(\theta) - \xtrue,\label{eq:hess-ip}
	\end{equation}
	i.e., a linear system whose coefficient matrix is the lower level Hessian. 
	Thus, a gradient-based method applied to \eqref{eq:ip-bl} involves the solution of a sequence of linear systems of the form \eqref{eq:hess-ip}, which stabilize as the iterations for solving \eqref{eq:ip-bl} proceed. 
	In applications such as inverse problems, each linear system \eqref{eq:hess-ip} is typically large scale and thus solved by iterative methods such as Krylov subspace methods. 
	Beyond warm-starting, a way to exploit similarity in a sequence of linear systems is to employ
	recycling Krylov subspace methods \cite{parks2006rks,soodhalter2020ssr,wang2007lto}, whereby 
	the solution space of one linear system is `recycled' in the solve of the next linear system.

	The remaining part of this paper is organized as follows. In Section~\ref{sec:background} we provide some background material on gradient-based methods for the solution of bilevel optimization problems, 
	as well as recycling Krylov methods for a sequence of linear systems. In Section~\ref{sec:recycle_what} we first briefly summarize existing strategies to determine  the subspace to be recycled, which are primarily based on 
	Ritz vectors, and then introduce Ritz generalized singular vectors that lead to a novel recycling strategy acknowledging that, in the bilevel optimization setting, solutions of the linear systems are employed in the construction of hypergradients. In Section~\ref{sec:stopcrit} we show how our proposed recycling strategy allows for a stopping criterion based on approximations of the hypergradient error, as opposed to existing stopping criteria based on the Hessian system residual in some norm \cite{gao2024ldl, pedregosa2016hoa,wang2007lto}. 
	Section~\ref{sec:numerics-krylov} reports the results of extensive numerical experiments performed on a couple of test problems, where we show the benefit of using recycling Krylov methods, and compare existing and our new recycling strategy and stopping criterion. Section \ref{sec:rec-conclusion} presents some concluding remarks. 
	
	
	\section{Background}\label{sec:background}
	We first show how a gradient-based method for bilevel optimization problems requires the solution of both the lower level problem and a linear system for the computation of the hypergradient at each iteration.
	We then discuss how the sequence of linear systems can be efficiently solved via  recycling Krylov subspace methods, with a particular emphasis on recycling MINRES, whereby a subspace from one linear system's approximation subspace is used in the solution of the next linear system.
	
	\subsection{Gradient-based solvers for bilevel optimization}
	\label{subsec:gb-methods}
	We consider the following class of bilevel problems:
	\begin{subequations}\label{eq:bl2}
		\begin{equation}
			\min_{\theta\in\real^p} \left\{\ulfun(\theta) = \ulfunx (\recon(\theta))\right\} \label{eq:ul2}
		\end{equation}
		\begin{equation}
			\text{subject to } \recon(\theta) = \argmin_{x\in\real^n} \llfun(x,\theta), \label{eq:ll2}
		\end{equation}
	\end{subequations}
	where $\ulfunx:\real^n\to\real$ is  differentiable, and $\llfun:\real^n\times\real^p\to\real$ is continuously differentiable with respect to $\theta,$ and twice continuously differentiable with respect to $x.$
	Additionally, we assume $\llfun$ is strongly convex with respect to $x,$ so that the lower level problem \eqref{eq:ll2} always admits a unique solution. 
	While  one could include a regularization term for $\theta$ in the the upper level cost $\ulfun,$ to simplify notation we do not consider this as it is straightforward to derive analogous results; see \cite{crockett2022bmi} for details.
	
	Application of the chain rule yields the gradient of the upper level cost function with respect to the hyperparameters, aptly named the hypergradient,
	\begin{equation}
		\nabla \ulfun(\theta) = \left(\deriv _\theta \recon(\theta) \right)^T \nabla \ulfunx(\recon(\theta)),\label{eq:ul-chain-rule}
	\end{equation}
	where $\deriv_\theta \recon(\theta)\in\real^{p\times n}$ is the Jacobian of $\recon(\theta)$ with respect to $\theta.$
	We typically select an upper level loss function such that $\nabla\ell$ is easy to compute, e.g., $\ell(x) = \tfrac{1}{2} \norm{x-\xtrue}^2_2$ yields $\nabla \ell(\recon(\theta)) = \recon(\theta) - \xtrue.$
	An issue with \eqref{eq:ul-chain-rule} is the application of the Jacobian $\deriv_\theta \recon(\theta).$ 
	Indeed, even with existence and uniqueness of $\recon(\theta),$ we often lack an analytic solution for which we can compute the corresponding Jacobian.
	However, by the assumptions on $\llfun,$  an alternative form of the Jacobian  motivated by the implicit function  theorem~\cite{crockett2022bmi,deoliveira2013iif,ehrhardt2023aih,zucchet2022bbo} can be employed. 
	Namely,
	\begin{equation}
		\deriv_\theta \hat x(\theta) = - \left[\deriv ^2 _{xx} \llfun(\hat x(\theta),\theta) \right]^{-1}  \deriv_{\theta x }^2\llfun(\hat x(\theta),\theta)\label{eq:implicit-update-context}
	\end{equation}
	where  $\deriv_{\theta x }^2\llfun := \deriv_{\theta }\nabla_x\llfun\in\real^{n\times p}$ and $\deriv^2_{xx}\llfun\in\real^{n\times n}$ is the Hessian of the lower level cost function which, by the assumptions on $\llfun,$ is SPD.
	Inserting \eqref{eq:implicit-update-context} into the hypergradient \eqref{eq:ul-chain-rule} yields 
	\begin{equation}
		\nabla \ulfun(\theta) = -\deriv_{\theta x }^2\llfun(\hat x(\theta), \theta)^T  \left[\deriv ^2 _{xx} \llfun(\hat x(\theta) , \theta) \right]^{-1}  \nabla \ulfunx(\recon(\theta)).\label{eq:hg-full}
	\end{equation}
	Regarding the computational feasibility of \eqref{eq:hg-full}, for applications such as imaging problems (whereby $x$ is a vectorized image) the Hessian is very large and, in general, lacks any structure.
	Indeed, we may only have access to a function that implicitly applies the action of the Hessian onto a vector.
	For these reasons the Hessian is not explicitly inverted in practice, and instead  an iterative method is employed to compute  the solution of the linear system
	\begin{equation}
		\deriv_{xx}^2\llfun(\reconp,\theta) \,w_\star = \nabla \ulfunx(\reconp)\label{eq:hess_sys}
	\end{equation}
	so that the hypergradient \eqref{eq:hg-full} is given by
	\begin{equation}
		\nabla\ulfun(\theta) = - (\deriv_{x\theta}^2\llfun(\reconp,\theta))^T w_\star.
		\label{eq:hypgrad}
	\end{equation}
	A class of methods that can be employed to solve \eqref{eq:hess_sys} are Krylov methods \cite{gazzola2020kmi, saad2003ims},  wherein an approximation subspace for the solution is generated via matrix-vector products with the Hessian, subject to constraints on the associated residual. A popular choice that exploits the SPD nature of the Hessian  is conjugate gradient~\cite{ crockett2022bmi,hestenes1952mcg,ji2021boc, pedregosa2016hoa}.

	At iteration $i$ of a gradient-based solver for \eqref{eq:bl2}, having the current solution $\theta^{(i)}$ and hypergradient \eqref{eq:hg-full} evaluated in $\theta^{(i)}$ (acting as a search direction $d^{(i)}$), an appropriate steplength can be determined via backtracking linesearch with the Armijo stopping rule. Namely, we find the smallest $j_i\in\mathbb{N}$ such that
	\begin{equation}
		\ulfun(\theta^{(i)} - \beta \rho^{j_i} d ^{(i)}) \leq \ulfun(\theta^{(i)}) - \eta \rho^{i_k} \nabla \ulfun (\theta ^{(i)})^T d^{(i)},\label{eq:armijo}
	\end{equation}
	where $\beta$ is an initial stepsize and $0< \eta, \rho <1$ and we eventually set $\theta^{(i+1)} = \theta^{(i)} - \beta \rho^{j_i} d ^{(i)}$. 
	Since evaluations of the upper level cost function $\ulfun$ are inherently expensive (each requiring a solution to the lower level problem \eqref{eq:ll2}), one must be sparing in how often a call is made when solving \eqref{eq:bl2}.
	Specifically, since the left hand side of \eqref{eq:armijo} is precisely $\ulfun(\theta^{(i+1)})$, it can be reused in the backtracking linesearch at iteration $i+1$ of the upper level solver, avoiding an unnecessary re-evaluation. Each intermediate evaluation of the left hand side in \eqref{eq:armijo} requires one lower level solve but, to reduce the computational burden of this task, we can employ warm restarting, wherein we initialise each lower level numerical solver with the solution of the previous solve. 
	The 
	reconstruction $\recon(\theta^{(i+1)})$ 
	can also be utilized in the construction of the hypergradient at iteration $i+1$ of the upper level solver. 
	
	An outline of backtracking linesearch that employs these efficiencies is provided in Algorithm~\ref{alg:backtrack}, and the pseudocode for gradient descent applied to \eqref{eq:bl2} is given in Algorithm~\ref{alg:hoag}.
	It should be noted that the hypergradient formula \eqref{eq:hg-full} we employ requires exact solutions to both the lower level problem and Hessian system. 
	Since we compute only numerical approximations to these subproblems, we can only estimate \eqref{eq:hg-full}.
	While a rigorous analysis of this aspect is beyond the scope of this paper, recent work has developed and analysed methods which acknowledge the inaccuracy of hypergradients \cite{ehrhardt2021ido,ehrhardt2023aih,pedregosa2016hoa,salehi2024aif}.

	\begin{algorithm}
		\caption{Backtracking linesearch for \eqref{eq:bl2}}\label{alg:backtrack}
		\begin{algorithmic}[1]
			\STATE Current solution $\theta,$ initial stepsize $\beta,$   search direction $d,$ factors $0<\mu,\rho<1$,
			
			\STATE If not already computed, determine $\tilde x_0 = \recon(\theta),$ $c_1 = \ulfun(\theta),$ and $c_2 = \nabla \ulfun(\theta)^Td$ 
			
			\FOR{$j=1,2,...$}
			\STATE Specify proposed update $\theta_j = \theta + \beta \rho^{j-1} d$
			\STATE Determine $\tilde x_j = \recon(\theta_j)$ using $\tilde x_{j-1}$ as an initial guess
			\STATE Evaluate $c_0 = \ulfun(\theta_j)$ using $\tilde x_j$
			\IF{Armijo condition $c_0 < c_1 + \eta \rho^{j-1} c_2$ is satisfied}
			\STATE Terminate, return $\theta_j,$ $\tilde x_j,$ next initial stepsize $\tilde\beta = 2 \beta\rho^{j-1},$ and $c_0$ 
			\ENDIF   
			\ENDFOR
		\end{algorithmic}
	\end{algorithm}

	\begin{algorithm}
		\caption{Gradient descent with backtracking for \eqref{eq:bl2}}\label{alg:hoag}
		\begin{algorithmic}[1]
			\STATE Initial guess $\theta^{(0)}\in\real^p,$ stopping gradient norm tolerance $\epsilon>0,$ 
			\STATE Determine $\recon(\theta^{(0)})$

			\FOR{$i=0,1,2,...$}
			\STATE Solve the linear system $\derivtwo_{xx} \llfun(\recon(\theta^{(i)}),\theta^{(i)}) \,w^{(i)}_\star = \nabla\ulfunx(\recon(\theta^{(i)}))$ \label{eq:code-hess}
			\STATE Construct search direction  $d^{(i)} = \derivtwo_{\theta x} \llfun(\recon(\theta^{(i)}),\theta^{(i)})^T\, w ^{(i)}_\star$
			\IF{$\norm{d^{(i)}}_2<\epsilon$}
			\STATE Terminate, return $\theta^{i)}$
			\ENDIF
			\STATE  Determine $\theta^{(i+1)}$ and $\recon(\theta^{(i+1)})$ via backtracking linesearch as in Algorithm~\ref{alg:backtrack}
			\ENDFOR
		\end{algorithmic}
	\end{algorithm}

	\subsection{(Recycling) Krylov methods} \label{subsec:krylov}
	The caveat of applying the solver described in the previous subsection is that, in order to form each hypergradient \eqref{eq:hypgrad}, one must determine the solution to both a sequence of lower level problems \eqref{eq:ll2} and Hessian system \eqref{eq:hess_sys}. 
	However, as the iterations of the gradient-based optimization scheme proceed, one would expect the computed hyperparameters to stabilise; in particular the associated sequences of solutions to the lower level and  Hessian systems would also stabilise. 
	
	We focus our attention on the sequence of linear systems that are solved. 
	To simplify the notation with respect to \eqref{eq:hess_sys}, we denote the Hessian system encountered at the $i$th iteration of a gradient-based method for solving \eqref{eq:bl2}, that is, line \ref{eq:code-hess} of Algorithm~\eqref{alg:hoag}, as
	\begin{equation}
		H^{(i)} w^{(i)}_\star = g^{(i)} \label{eq:hess_seq}
	\end{equation}
	where $H^{(i)}:= \deriv_{xx}\llfun(\reconpk,\theta^{(i)})\in\real^{n\times n}$ and $g^{(i)}:=\grad_x\ulfunx(\reconpk)\in\real^n.$
	In practice, as explained in Section~\ref{subsec:gb-methods}, we can only compute an approximate solution of  \eqref{eq:hess_seq}, denoted by $w^{(i)},$ via an iterative solver. 
	One way to yield computational speedup in solving the sequence of linear systems \eqref{eq:hess_seq} is via warm restarting, wherein the iterative solver for the $i$th problem is initialised with the determined solution of the $(i-1)$th system, $w^{(i-1)}.$ 
	However, if the iterative solver for \eqref{eq:hess_seq} is a subspace projection method (e.g., a Krylov method), 
	similarity of the linear systems can be exploited further via recycling \cite{parks2006rks,soodhalter2020ssr,wang2007lto}, a technique in which  information from the  solution subspace of the $(i-1)$th system solve is used to kickstart the solve for the $i$th system.
	Specifically, the idea of recycling is the following: 
	\begin{itemize}
		\item Approximate the solution of the first  Hessian system of \eqref{eq:hess_seq} with 
		$$w^{(0)} = w_0 + V^{(0)} y^{(0)}$$
		for initial guess $w_0\in\real^n$ and coefficients $y^{(0)}\in\real^{k_0}$ associated with a basis $V^{(0)}\in\real^{n\times k_0}$  of the generated Krylov subspace.
		\item From $V^{(0)}$ extract `useful' vectors $U^{(1)}\in\real^{n\times s_1}$ and demand that the solution of the second system in \eqref{eq:hess_seq} is of the form 
		$$w^{(1)} = w^{(0)} + V^{(1)} y^{(1)} + U^{(1)} z^{(1)},$$
		where $y^{(1)}\in\real^{k_1}, z^{(1)}\in\real^{s_1},$  and the columns of  $V^{(1)}\in\real^{n\times k_1}$ are a basis for the generated Krylov subspace.
		\item For subsequent linear systems $i=2,3,\ldots,$ extract `useful' vectors $U^{(i)}\in\real^{n\times s_i}$ from  $[V^{(i-1)}, \;U^{(i-1)}] $ and similarly demand that approximate solutions are of the form
		$$w^{(i)} = w^{(i-1)} + V^{(i)} y^{(i)} + U^{(i)} z^{(i)},$$
		where  $y^{(i)}\in\real^{k_i}, z^{(i)}\in\real^{s_i},$ and the columns of  $V^{(i)}\in\real^{n\times k_i}$ are a basis for the generated Krylov subspace.
	\end{itemize}
	
	Each $U^{(i)}$ is called a recycle space, since information from the previous solution space is `recycled' for the current linear system solve.
	It is crucial that, given $U^{(i)},$ one can generate suitable basis vectors $V^{(i)}$ and specify coefficients $y^{(i)}$ and $z^{(i)},$ such that $w^{(i)}$ is a meaningful approximation of $w_\star^{(i)}.$
	Current recycling methods primarily consist of Krylov methods  that have been carefully adapted to account for a recycle space \cite{bolten2022ksr,wang2007lto}. 
	In this work we focus on the RMINRES \cite{wang2007lto}, a recycling variant of MINRES, a Krylov subspace method in-which the symmetry of the Hessian is exploited to minimize the residual norm.

	\subsubsection{MINRES}
	
	We now focus on how a given linear system in \eqref{eq:hess_seq} can be solved.
	To further simplify the notation, we suppress the superscript in \eqref{eq:hess_seq} and consider the linear system 
	$Hw_\star=g,$ with initial guess $w_0.$ 
	One way to solve the linear system is via a Krylov subspace method \cite{gazzola2020kmi,saad2003ims}, such as MINRES \cite{paige1975ssi}, which builds a Krylov subspace 
	$\kry_k(H,r_0) := \mathrm{span}\{r_0,Hr_0,\dots, H^{k-1}r_0\},$ where $r_0=g - Hw_0$ is the initial residual, and computes a solution that is optimal over $\kry_k(H,r_0).$ 
	
	Specifically, MINRES determines a correction term $\xi_k\in\kry_k(H,r_0)$ such that $w_k = w_0 + \xi_k$ minimizes $\norm{r_k}_2.$
	To do this, symmetry of $H$ is exploited to efficiently generate  an orthonormal basis of $\kry_k(H,r_0)$ via the Lanczos relation \cite{saad2003ims}:
	\begin{align}
		v_1 = r_0, \quad \beta_1 = \norm{v_1}_2\nonumber
		\\
		\beta_{j+1}v_{j+1} = Hv_j - \alpha_jv_j - \beta_j v_{j-1}, \label{eq:arnoldi-3-term-rec}
	\end{align}
	where	$\alpha_j := v_j^T H v_j$ and $\beta_{j+1} := \norm{v_{j+1}}_2$ for $j=1,\cdots,k.$ 
	The Lanczos relation can be written compactly as
	\begin{equation}
		HV_k = V_{k+1} \bar T_k,\label{eq:arnoldi-matrix-id}
	\end{equation}
	where  $V_k = [v_1\;\cdots \; v_k ]\in\real^{n\times k},$ $V_{k+1}=   [V_k\; v_{k+1} ]\in\real^{n\times(k+1)}$ have orthonormal columns such that $\text{range}(V_k)=\kry_k(H,r_0)$,  the leading minor $T_k$ of order $k$ of $\bar T_k \in\real^{(k+1)\times k}$ is symmetric and tridiagonal with $\alpha_1,\cdots ,\alpha_k$ along the diagonal and $\beta_2,\cdots \beta_{k}$ on the off-diagonals, and $\bar T_k =[T_k^T, \beta_{k+1}e_{k}]^T$, with $e_k\in\real^{k}$ being the $k$th canonical basis vector. Here and in the following we assume that the Lanczos algorithm does not break down, i.e.,the  that relation \eqref{eq:arnoldi-matrix-id} can always be written. 
	Then,  exploiting the Lanczos relation \eqref{eq:arnoldi-matrix-id}, $\xi_k = V_ky_k$ is determined by solving
	\begin{align*}
		y_k = \argmin_{y\in\real^k} \norm{HV_ky - r_0}^2_2
		= \argmin_{y\in\real^k} \norm{V_{k+1} \bar T_k y - \beta_1 v_1}^2_2
		= \argmin_{y\in\real^k} \norm{\bar T_k y - \beta_1 e_1}^2_2
	\end{align*}
	where $e_1\in\real^{k+1}$ is the first canonical basis vector. 
	If the residual norm is sufficiently small the algorithm can be terminated, otherwise a new basis vector can be efficiently constructed by the 3-term recurrence \eqref{eq:arnoldi-3-term-rec} and the associated residual norm determined. In particular, MINRES reduces the solution of the original linear system to the solution of a linear least squares problem with tridiagonal coefficient matrix of order $k$, which can be efficiently performed; indeed, the Lanczos relation \eqref{eq:arnoldi-3-term-rec} allows the QR decomposition of $\bar T_k$ to be efficiently updated to yield the QR decomposition of $\bar T_{k+1}$; see \cite{paige1975ssi} for details.

	\subsubsection{RMINRES}\label{subsec:arnoldi_recycle}
	RMINRES \cite{wang2007lto} is a recycling variant of MINRES that finds a correction term of the form
	\begin{equation}
		\xi_k = V_k y_k + U z_k \label{eq:correction_term}
	\end{equation}
	such that 
	\begin{equation}
		(y_k,z_k) =	\argmin_{y\in\real^k,\,z\in\real^s} \norm{r_0 - H\begin{bmatrix} V_k & U \end{bmatrix}  \begin{bmatrix}y\\z \end{bmatrix}}^2_2,\label{eq:coupled}
	\end{equation}
	where $U\in\real^{n\times s}$ is a pre-specified recycle space and the columns of $V_k\in\real^{n\times k}$ form an orthonormal basis of a Krylov subspace which is iteratively constructed.
	We remind ourselves that, although we suppress the notation in this subsection, in the framework of this paper every variable in \eqref{eq:correction_term} - including the recycle space $U$ - depends on the iteration $(i)$ of the gradient-based algorithm employed to solve the bilevel problem.

	Since the correction term \eqref{eq:correction_term} is partially represented by the columns of $U$ by construction, it would be a waste of computational resources to allow the  generated columns of $V_k$  to overlap with information already captured by $U.$  
	To this end, RMINRES ensures each Lanczos vector $v_j$ is orthogonal to $C:=HU,$ where $U$ is chosen such that $C^TC=I,$  which results in the augmented the Lanczos relation  
	\begin{equation}
		(I-CC^T)HV_k = V_{k+1} \bar T_k,\label{eq:arnoldi_aug_matrix}
	\end{equation}
	where $\bar T_k$ can be shown to be tridiagonal \cite[Section~5.2]{soodhalter2020ssr}.
	This precise choice of orthogonalisation allows
	\eqref{eq:coupled} to be decoupled into
	\begin{equation}
		y_k = \argmin_{y\in\real^k} \norm{\beta_1 e_1 - \bar T_k y}_2^2 \label{eq:y_solve}
	\end{equation}
	\begin{equation}
		z_k =C^Tr_0 - C^THV_ky_k.\label{eq:zk_answer}
	\end{equation}
	In particular, \eqref{eq:y_solve} has the same tridiagonal structure as the problem considered by standard MINRES, and thus can be solved in a similarly efficient manner with minor modifications to account for the part of the solution in the recycle space  with coefficients \eqref{eq:zk_answer}; 
	see \cite{wang2007lto} for details.
	
	We quickly discuss how one can ensure that $C=HU$ satisfies $C^TC=I.$ 
	It is important to remind ourselves that we only care about the column space (range) of the recycle matrix, rather than the columns themselves.
	With this in mind, suppose we have determined that 
	$\mathrm{range}(\tilde U)$
	for some $\tilde U\in\real^{n\times s}$  is a good recycle space. 
	Then, consider the reduced QR decomposition $
	H\tilde U = CR$
	where $C\in\real^{n\times s}$ has orthogonal columns and $R\in\real^{s\times s}$ is upper triangular. 
	We see that  $U:= \tilde U R^{-1}$  clearly satisfies the condition that we want while losing no information since, by the invertability of $R,$ $\mathrm{range}(U) = \mathrm{range}(\tilde U).$
	In the following discussion we assume that this QR decomposition is implicitly performed, and make no distinction between $U$ and $\tilde U,$ referring to both as the recycle space.
	

	\section{Defining and computing the recycle space for bilevel optimization}\label{sec:recycle_what}

	We desire a recycling strategy such that the sequence of Hessian systems \eqref{eq:hess_seq} can be quickly solved, and such that the associated hypergradients  \eqref{eq:hypgrad} are accurate.
	In Section~\ref{subsec:existing-recycle-strats} we cover  existing recycling strategies \cite{garcia2024irm,morgan2000irg,soodhalter2020ssr,wang2007lto,yetkin2023rnk}, which use eigenvector approximations called Ritz vectors.
	While such strategies can be employed to reduce the cost associated with solving the sequence of Hessian systems \eqref{eq:hess_seq}, they only consider the sequence in isolation of how the solutions are utilized.
	Indeed, in the hypergradient \eqref{eq:hypgrad}, solutions of the linear systems are premultiplied by the adjoint Jacobian  of the gradient of the lower level cost function.
	Thus, investment in an accurate or computationally cheap $w^{(i)}$ is only a proxy for an accurate or computationally cheap hypergradient.
	Motivated by this, we propose in Section~\ref{subsec:recycle-bilevel} a novel  recycling strategy which account for this pre-multiplication, and thus acknowledge the bilevel setting in which the linear systems arise from.
	
	\subsection{Existing selection criteria}\label{subsec:existing-recycle-strats}
	Existing works that employ recycling strategies are primarily motivated by two settings: a sequence of slowly changing linear systems are solved; or a single linear system is solved but, due to computational memory restrictions, the numerical solver must be restarted mid-solve~\cite{morgan2000irg}.
	As we are interested in solving \eqref{eq:hess_seq}, we tailor discussion to the former.
	With this in mind, we interpret $H$ as the Hessian of the system we are about to solve, i.e., $H=H^{(i)};$  we include a discussion on the implications of this choice later in the section.
	
	Although the two settings are different, existing recycling strategies that exploit approximate invariant subspaces are motivated by approximating the eigendecomposition 
	\begin{equation}
		H = \widehat Q\widehat \Lambda \widehat Q ^T, \label{eq:full-eig-decomp}
	\end{equation}
	where $\widehat Q\in\real^{n\times n}$ is orthonormal and $\widehat \Lambda\in\real^{n\times n}$ is a diagonal matrix of the eigenvalues.
	Indeed, if one has access to the full eigendecomposition \eqref{eq:full-eig-decomp} then the correction term is directly $\xi = \widehat Q\widehat \Lambda^{-1} \widehat Q^T r_0.$
	Due to the large-scale nature of the system and, in general, lack of structure of the Hessian, it is not feasible to compute the eigendecomposition for every Hessian system, hence an iterative method is employed. 
	Even so, since $\xi$ is in the column space of $\widehat Q,$ choosing $U$ to mimic $\widehat Q$ in some capacity is  a natural choice. 
	
	One idea is to consider the eigendecomposition of $H$ restricted to the range of some matrix $W\in\real^{n\times t},$ $n\gg t,$ for example, a basis of the Krylov subspace generated by the previous Hessian system solve.
	To do this, one can compute the eigendecomposition 
	\begin{equation}
		W^THW = Q\Lambda Q^T \in\real^{t\times t} \label{eq:ritz-matrix-decomp}
	\end{equation}
	and consider a subspace in the range of $ WQ.$
	The columns of $WQ$ are called Ritz vectors, and the associated eigenvalues (entires of the diagonal matrix $\Lambda$) are called the Ritz values. 
	Notice that should $n=t,$ we recover the full eigendecomposition
	$
	H = WQ\Lambda (WQ)^T.
	$
	Ritz values tend to approximate the extremal eigenvalues of $H$ well, but can give poor approximations to interior eigenvalues 	\cite{nakatsukasa2017asv}.
	For example, eigenvalues can be written as a min-max characterisation of the Rayleigh quotient ~\cite[Theorem 10.2.1]{parlett1998sep} and one of the Ritz values is the largest possible Rayleigh quotient over the given subspace ~\cite[Theorem 11.4.1]{parlett1998sep}.
	Other characterizations of how Ritz vectors/values can be considered optimal approximations of eigenvectors/values over the subspace spanned by $W$ are given in 
	\cite[Chapter~8]{genehgolub2013mc}.
	
	It turns out that the Lanczos algorithm for generating Krylov basis vectors is the Rayleigh-Ritz procedure applied to a Krylov subspace \cite[Chapter 13]{parlett1998sep}, and so the generated Krylov subspace inherently approximates exterior eigenvectors. 
	For this reason, recycling approximations of interior eigenvectors can propagate information least readily available from the Lanczos algorithm.
	An approximation of interior eigenvalues are Harmonic Ritz values, which are the reciprocal of (ordinary) Ritz values of $H^{-1}$ with respect to $ W.$ 
	Harmonic Ritz values are named as such due to the interpretation of them as a weighted harmonic mean of $H$'s eigenvalues, in contrast to Ritz values which can be viewed as a weighted arithmetic mean of $H$'s eigenvalues; see \cite{paige1995ase} for details.
	In practice,  harmonic Ritz vectors are not determined by forming $H^{-1}$ and computing the associated Ritz vectors, but instead an equivalent generalized eigenvalue problem
	\begin{equation}
		(HW)^THW\rho = \tilde \theta (HW)^T W \rho\label{eq:harmonic-ritz}
	\end{equation}
	is solved, and harmonic Ritz vectors are given by $W\rho$ \cite{wang2007lto}.

	Besides recycling (harmonic) Ritz vectors, other recycling strategies are possible. For example, the authors of \cite{kilmer2006rsi} adopt recycle spaces obtained by combining subspaces of past
	approximate solutions with approximate invariant subspaces. A recycling strategy inspired by goal-oriented proper orthogonal decomposition from model reduction is proposed in \cite{Tuminaro}. The authors of \cite{AlDaas} use a recycle space defined by approximate right singular vectors of the coefficient matrix. There are also a number of recycling or subspace augmentation strategies explored in the setting of inverse problems, where the main principle underlying the latter is to enrich the approximation subspace for the solution by vectors that may enhance the quality of the reconstruction. For instance, in  \cite{jiang2021hpm, RMMGKSp} the recycle space is selected according to specific compression techniques applied to either the projected coefficient matrix or the approximation subspace $\text{range}(V_k)$ for the solution; among the latter, choosing the columns of $V_k$ that correspond to the coefficients in $y_k$ of largest magnitude allows to retain the most important columns of $V_k$ that represent the solution. An effective augmented GMRES version was first proposed in \cite{baglama2007augmented} and further investigated, among other solvers and even in a continuous setting, in \cite{ramlau2021subspace}. Similarly, there is an extensive literature about the so-called generalized Krylov subspace methods, where the standard Krylov subspace for the solution (associated to GMRES or CGLS) either appears alongside additional vectors (like in \cite{reichel2012tikhonov}) or is modified by  multiplication with iteration-dependent matrices (like in \cite{buccini2023}), which are linked to Tikhonov-regularized problems.

	With regards to solving a sequence of linear systems, the $i$th recycle space $U^{(i)},$ which is employed to solve the $i$th linear system $H^{(i)}  w^{(i)} = g^{(i)},$ is determined using  $W^{(i)}=[V^{(i-1)}\;U^{(i-1)}],$ that is, information from the $(i-1)th$ linear system $H^{(i-1)}w^{(i-1)}=g^{(i-1)}.$
	A subtle detail is whether $H^{(i-1)}$ or $H^{(i)}$ should be used alongside  $W^{(i)}$ to determine $U^{(i)}.$
	For example, to compute Ritz vectors one could consider the eigendecomposition of either
	\begin{align}&(W^{(i)})^T  H^{(i)} W^{(i)}.\label{eq:subtle-new}
		\\ \text{or}\quad&  (W^{(i)}) ^T H^{(i-1)} W^{(i)}\label{eq:subtle-old}
	\end{align}
	Since we desire a recycle space that is appropriate for solving the $i$th linear system,  \eqref{eq:subtle-new} is a natural consideration.
	However, \eqref{eq:subtle-new} requires that the entire basis of the Krylov subspace is retained between linear system solves, even if the Krylov method itself does not demand this.
	For this reason, some existing work \cite{wang2007lto} advocate \eqref{eq:subtle-old}, as then one can not only iteratively update the recycle space to be used in the next linear system during the current linear system solve, but also exploit known  matrix structure and decompositions associated with $H^{(i-1)} W^{(i)},$ such as the Lanczos relation \eqref{eq:arnoldi_aug_matrix}.
	Unless otherwise stated, we employ \eqref{eq:subtle-new} and leave a thorough  analysis of using \eqref{eq:subtle-old} as future work.

	\subsection{Acknowledging the bilevel setting}\label{subsec:recycle-bilevel}
	In the bilevel setting, we are interested in the hypergradient
	\begin{equation*}
		\nabla \ulfun(\theta) = \jacadj H^{-1} g,
	\end{equation*} 
	where $\jacadj:=  - (\deriv_{\theta x}\llfun(\reconp,\theta))^T\in\real^{p\times n}.$
	In practice we first approximate $w_\star = H^{-1}g$ and then perform a premultiplication by $\jacadj.$ 
	The recycling strategies outlined in Section~\ref{subsec:existing-recycle-strats} use vectors related to projected approximations of $H^{-1}$ and so are oblivious to the premultiplication applied to the determined solution.
	To acknowledge the bilevel setting, 
	we propose to approximate $\jacadj H^{-1}$ over the column-space of a full-rank matrix $W\in\real^{n\times t}$ by using the  Generalized Singular Value Decomposition (GSVD) \cite{bjorck1990lsm,genehgolub2013mc,vanloan1976gsv},
	which provides a simultaneous as defined in the following statement. 
	
	\begin{theorem}[GSVD]\label{thm:gsvd}
		Let  $\tilde  H\in\real^{t\times t}$ be invertible and  $\tilde \jacadj\in\real^{p\times t}$ 	with $p\geq t.$
		There exist orthogonal $V_{\tilde \jacadj} \in\real^{p\times p}$ and $V_{\tilde H}\in\real^{t\times t}$ and invertible $X\in\real^{t\times t}$ such that	
		\begin{align*}
			V_{\tilde \jacadj}^T \tilde \jacadj X = D_{\tilde \jacadj} = 
			\begin{matrix}
				\begin{bmatrix}
					\diag(\alpha_{1},\ldots, \alpha_t )  \\
					0 
				\end{bmatrix}
			\end{matrix} \in\real^{p\times t},
			\\
			V_{\tilde H}^T \tilde H X  = D_{\tilde H}
			= \diag(\beta_{1},\ldots, \beta_t ) \in\real^{t\times t},
		\end{align*}	
		where $	0\leq \alpha_{1} \leq \ldots \leq \alpha_t< 1$ and $1\geq \beta_{1} \geq \ldots \geq \beta_t> 0$ are such that
		\begin{align*}
			\alpha_i^2 + \beta_i^2 = 1,\quad i=1,\cdots ,t.
		\end{align*}

		The quantities $ \mu_i = \alpha_i / \beta_i,\; i=1,\cdots t$ are called the generalized singular values and are such that $0\leq \mu_1\leq \mu_2\leq \ldots \leq \mu_t.$ 
		We refer to $X$ as the right singular vectors, and call $V_{\tilde \jacadj}$ and $V_{\tilde H}$ the left singular vectors of $\tilde \jacadj$ and $\tilde H$ respectively.
		
	\end{theorem}
	The GSVD reverts to the SVD should $\tilde \jacadj=I\in\real^{t\times t}$, since we can set $X=V_{\tilde\jacadj}$  to recover $V_{\tilde H} ^T\tilde H V_{\tilde \jacadj} = D_{\tilde H},$ i.e., the SVD of $\tilde H.$
	While we assume $p \geq t$ both in Theorem~\ref{thm:gsvd} and in the following,  analogous results  follow for the case $p<t.$
	Similarly, the GSVD can be stated for singular $\tilde H,$ but this will not be the case for our setting.
	
	By considering the GSVD of the pair $(\jacadj W, W^THW)$ 
	\begin{equation}
		\jacadj W = V_\jacadj D_\jacadj X^{-1},\quad W^THW=V_H D_H X^{-1}\,, \label{eq:gsvd}
	\end{equation}
	it follows that
	\begin{equation}
		\jacadj W (W^THW)^{-1}  = V_\jacadj D_\jacadj D_H ^{-1} V_H^T \in\real^{n\times t} . \label{eq:gritz-decomp}
	\end{equation}
	Should $W$ be square and orthogonal, \eqref{eq:gritz-decomp} recovers the SVD of $\jacadj H^{-1} W$ which, by orthogonality of $W,$ is equivalent to the SVD of $\jacadj H^{-1}.$

	It follows from \eqref{eq:gsvd} that we have the diagonalisation $(WV_H)^TH WX = D_H.$
	With this in mind, we call the columns of $WV_H$ the left Ritz  generalized  singular vectors, the columns of  $WX$ the right  Ritz  generalized singular vectors, and the diagonal entries of $D_\jacadj D_H^{-1}$ the Ritz generalized singular values. 	
	Notice that, if we take $\jacadj W=I,$ then the Ritz generalized singular values simply are the reciprocals of the Ritz values of $H.$
	Unlike recycling using Ritz vectors, where there was only a single choice regarding what vectors  could be used as the recycle space, here we could consider:
	the left Ritz generalized singular vectors, $ WV_H;$
	the right Ritz generalized singular vectors, $WX;$
	a linear combination of both the  left and right Ritz generalized singular vectors, e.g. $\tfrac{1}{2}(WV_H + WX).$
	
	Then the recycle space $U\in\real^{n\times s}$ is given by the columns of the selected matrix corresponding to $s$ among the smallest and/or largest Ritz generalized singular values. 
	Thus, in practice, we need only compute the partial GSVD of  $(\jacadj W, W^THW)$ for $s$ components \cite{alvarruiz2024tjl,huang2022thj,jia2023cgc}, which finds matrices $\tilde V_\jacadj \in\real^{p\times s}, \tilde V_H \in\real^{t\times s}$ with orthogonal columns, diagonal matrices $\tilde D_\jacadj\in\real^{s\times s},$ $\tilde D_H\in\real^{s\times s},$ and full rank matrix  $\tilde X\in\real^{t\times s}$ such that
	\begin{align}
		\jacadj W  \tilde X = \tilde V_\jacadj \tilde D_\jacadj,
		\qquad 
		W^THW \tilde X = \tilde V_H \tilde D_H.\label{eq:partial-gsvd}
	\end{align}
	
	We conclude this section by emphasising again that building the recycle space upon the newly-introduced Ritz generalized singular vectors and values provides a heuristic way of simultaneously evaluating the action of $J$ on the subspace $\text{range}(W)$, which contains the approximate solution $w$ to the current Hessian system. Although not considered in this paper, there may be other strategies to assess the impact of premultiplying 
	$w$ by $J$. For instance, following the optimal truncation results presented for a (sequence of) linear system(s) in \cite{desturler1999tso, parks2006rks}, it may be possible to prove that retaining the $s$-dimensional subspace of $\text{range}(W)$ that better aligns (in terms of principal angles between subspaces \cite{genehgolub2013mc}) with $\text{range}(J^T)$ leads to an optimal approximation of $Jw$ within any $s$-dimensional subspace of $\text{range}(W)$. In practice, 
	assuming that computing the SVD of $J$ is feasible, and denoting by $V_J$ its right singular vector matrix, one may choose $U=WV_s$, where $V_J^TW=U_s\Sigma_sV_s^T$ is the rank-$s$ truncated SVD of $V_J^TW$. Ways of achieving this in a matrix-free setting will be the focus of future research.

	\section{Alternative stopping criterion for \eqref{eq:hess_seq}}\label{sec:stopcrit}
	Recall that the standard \linebreak[4] (R)MINRES algorithm is terminated when the residual satisfies
	\begin{equation*}
		\delta >\norm{Hw_k - g}_2  =: \norm{r_k}_2
	\end{equation*}
	for some tolerance $\delta>0.$
	Ideally, we would terminate the iterative solver once the associated hypergradient is sufficiently accurate, that is, once
	\begin{align}
		\delta &> \norm{\jacadj w_k - \jacadj w_\star}_2 
		= \norm{\jacadj H^{-1} r_k}_2.\label{eq:hg-error-stop}
	\end{align}
	In order to approximate \eqref{eq:hg-error-stop} notice that we can rewrite the second equation in \eqref{eq:partial-gsvd} as $(W^THW)^{-1}\tilde V_H=\tilde X\tilde D_H^{-1} .$
	With this in mind, consider the approximation
	\begin{align*}
		\jacadj H^{-1} 
		&\approx \jacadj W (W^THW)^{-1}W^T 
		\approx
		\jacadj W \tilde X \tilde X^\dagger (W^THW)^{-1} \tilde V_H \tilde V_H^T W^T 
		\\
		&= \jacadj W \tilde X \tilde X^\dagger \tilde X \tilde D_H^{-1} \tilde V_H^T W^T 
		=
		\tilde V_\jacadj \tilde D_\jacadj \tilde D_H^{-1} \tilde V_H^T W^T,\label{eq:gritz-decomp-lowrank}
	\end{align*}
	where the final equality follows by the definition of the pseudoinverse and \eqref{eq:partial-gsvd}.
	In the case that the full GSVD is employed and $W$ is square, we have equality throughout the above approximation of $\jacadj H^{-1}.$
	Thus, we approximate \eqref{eq:hg-error-stop} with
	\begin{align}
		\delta >  \norm{\tilde D_\jacadj \tilde D_H^{-1} \tilde V_H^T W^T r_k}_2
		\label{eq:hg-err-stop-lowrank}.
	\end{align}

	We want to emphasise that, although Ritz vectors can be connected to projected eigenvectors, we have not shown that an analogous statement can be made for Ritz generalized singular vectors and generalized singular vectors. 
	Furthermore, we cannot conclude that the hypergradient error approximation \eqref{eq:hg-err-stop-lowrank} is a projection of \eqref{eq:hg-error-stop}.
	
	We remark that to employ \eqref{eq:hg-err-stop-lowrank} we require an efficient way to determine the current residual vector $r_k.$ 
	While RMINRES exploits matrix decompositions to iteratively update the residual norm directly and avoid explicit computation of the residual vector, one can utilize said matrix decompositions encountered in RMINRES to efficiently determine $r_k$ with little overhead; see supplementary material.

	\section{Numerical experiments}\label{sec:numerics-krylov}
	We focus on a bilevel learning problem motivated by inverse problems, where the task is to learn optimal filters of the Fields of Experts (FoE) \cite{roth2009fe}  regularized problem 
	\begin{equation}
		\min_{x\in \real^n} \left\{\frac{1}{2}\norm{\op x - \ynoise}^2_2 + \frac{\epsilon}{2}\norm{x}^2_2   +   \sum_{i=1}^N e^{\theta_0^{(i)}} 
		\left\langle 1 , \phi(\tilde \theta^{(i)}*\,x) \right\rangle\right\},
		\label{eq:ll2ex} 
	\end{equation}
	where  $*$ denotes a convolution,  $1\in\real^n,$ and $\phi$ is an elementwise operation (e.g., $\phi(s)= s^2$ recovers the squared $\ell_2$-norm).
	The parameters to be learned are the exponents  $\theta_0^{(i)}\in\real$ defining the regularization parameters and the filters $\tilde \theta ^{(i)}\in\real^{\tilde p_i}$ for $i=1,\cdots,N,$ so that $\theta = [\theta_0^{(1)},(\tilde\theta^{(1)})^T,\ldots,\theta_0^{(N)},(\tilde\theta^{(N)})^T]^T\in\real^p.$
	To guarantee existence and uniqueness of solutions to \eqref{eq:ll2ex}, we fix $\epsilon=10^{-6}.$ 
	
	We consider $K$ training pairs $(\xtrue_{j}, y_{j})$, $j=1,\dots,K$ and focus on  the bilevel learning problem
	\begin{equation}
		\min_{\theta\in\real^p} \frac{1}{K}\sum_{j=1}^K \norm{\recon_{j}(\theta) - \xtrue_{j}}_2^2,\; \text{where each }\recon_{j}(\theta)\text{ solves \eqref{eq:ll2ex} with }y=y_{j}. \label{eq:ul2ex}
	\end{equation} 
	
	We investigate two applications: an inpainting problem of the MNIST dataset \cite{yannlecun2010mhd}, wherein the relatively small dimension of the inverse problem allows for a rigorous empirical exploration; a deconvolution problem of the BSDS300 dataset \cite{MartinFTM01}, to illustrate that the proposed recycling strategies remain effective for larger problems.

	Before describing each test in detail, we describe 
	the common computational framework employed for our numerical experiments. 
	We numerically solve the bilevel problem utilizing various recycling strategies, and our Python implementation   utilizing PyTorch  \cite{paszke2019pis} can be found  at \href{https://github.com/s-j-scott/bilevel-recycling}{github.com/s-j-scott/bilevel-recycling}. 
	To the best of our knowledge, unlike the partial eigendecomposition and generalized eigendecomposition, there is no readily available general-purpose implementation of the (partial) GSVD in Python.
	Thus, to employ our proposed recycling strategy, we 
	follow the GSVD construction as in \cite{paige1981gsv}, which involves the computation of two SVDs and a QR decomposition. 
	While (partial) matrix decompositions can in principle be performed in a matrix-free fashion (i.e., via a function that returns the respective matrix-vector product), the implementations we utilize require explicit access to the full matrix.
	Thus, we construct full matrix representations of $\jacadj W$ and $W^THW$ in PyTorch.
	We remark that, to compute Ritz vectors, which only concerns an eigendecomposition, we still compute a full matrix representation of $W^THW$, but then perform the partial eigendecomposition in Python via Scipy.

	If we were to solve the entire bilevel problem from scratch with each recycling strategy, a new sequence of Hessians would technically be generated for each strategy.
	Therefore, in order to have a fair comparison of the performance of the different recycling strategies 
	we require a single sequence of Hessian systems that can be solved using each strategy. 
	To do this, we solve the bilevel problem \eqref{eq:ul2ex} with $K=1$, without any recycling and tolerances to be specified.

Given a solution $w$ to the $i$th linear system, we want to measure the quality of the determined hypergradient (abbreviated to HG in the plots), $\jacadj^{(i)} w.$
To do this,  we compute high-quality reference hypergradients $\jacadj^{(i)} w_\star^{(i)},$ where $w_\star^{(i)}$ was obtained by solving the $i$th linear system without recycling and with a residual norm stopping criterion tolerance of $\delta =10^{-13}.$
We can then compute the relative error of the hypergradient associated with the $i$th iteration of the upper level solver, given by
${\|{\jacadj^{(i)} w_\star^{(i)} - \jacadj^{(i)} w\|}_2}/{\|{\jacadj^{(i)} w_\star^{(i)}\|}_2}.$

We are also interested in the potential impact the hypergradients determined via recycling may have on the parameters and associated upper level cost. 
Since we already solve a fixed sequence of Hessians associated to the bilevel problem without recycling, we implicitly already have a sequence of parameters $\theta^{(i)},$ one for each linear system.
To determine the influence of the hypergradient obtained via recycling, for the $i$th linear system and a given recycling strategy, we start with $\theta^{(i)}$ determined from the non-recycling run of the sequence of systems, employ recycling to solve the associated Hessian for $\hat w^{(i)},$ calculate the corresponding hypergradient $\jacadj^{(i)} \hat w^{(i)},$  and then use backtracking linesearch to determine the updated parameter $\hat  \theta^{(i+1)}$ and associated upper level cost. 
We remark that, in general, $\hat \theta^{(i+1)}$ is specific to the recycling strategy and can be slightly different to $\theta^{(i+1)}.$

To simplify the notation, we denote recycling strategies to be of the form \textit{VecType-Size(Side)} where
\begin{itemize}
	\item \textit{VecType} denotes the sort of vector quantity utilized: Eigenvectors (Eig), Ritz vectors (Ritz), harmonic Ritz vectors (HRitz), generalized singular vectors (GSVD), Ritz generalized singular vectors (RGen).
	\item \textit{Size} denotes whether we select vectors corresponding to the smallest (S), largest (L) or mixture of the smallest and largest (M) eigenvalues/(harmonic) Ritz values/ (Ritz) generalized singular values.
	If we consider a mixture, we select half the dimension of the recycle space number of vectors to be associated with the small quantities, and half with the large.
	\item When relevant, \textit{Side} denotes whether we select right (R), left (L) or a mixture of right and left (M) (Ritz) generalized singular vectors.
\end{itemize} 
For example, recycling using the Ritz vectors associated with the smallest Ritz values is denoted Ritz-S, and recycling using the right Ritz generalized singular vectors associated with the largest Ritz generalized singular values is denoted RGen-L(R). Table~\ref{tab:acro} summarizes all the acronyms used in this section along with their meanings.

\begin{table}
	\caption{Explanation of the acronyms used to signify different recycling strategies.}\label{tab:acro}
	\centering
	\begin{tabular}{cclc}\hline
		Acronym &Side &  VecType & Size\\\hline
		Eig-S &- & Eigenvectors & Small \\
		GSVD-L(R) &Right & Generalized singular vectors & Large \\
		Ritz-S & - & Ritz vectors & Small \\
		Ritz-L & - & Ritz vectors & Large \\
		Ritz-M & - & Ritz vectors & Small \& large \\
		HRitz-S & - & Harmonic Ritz vectors & Small \\
		HRitz-L & - & Harmonic Ritz vectors & Large \\
		HRitz-M & - & Harmonic Ritz vectors & Small \& large \\
		RGen-S(R) & Right & Ritz generalized singular vectors &  Small\\
		RGen-L(R) & Right & Ritz generalized singular vectors &  Large\\
		RGen-M(R) & Right & Ritz generalized singular vectors &  Small \& large\\
		RGen-S(L) & Left & Ritz generalized singular vectors &  Small\\
		RGen-L(L) & Left & Ritz generalized singular vectors &  Large\\
		RGen-M(L) & Left & Ritz generalized singular vectors &  Small \& large\\
		RGen-S(M) & Left \& right & Ritz generalized singular vectors &  Small\\
		RGen-L(M) & Left \& right & Ritz generalized singular vectors &  Large\\
		RGen-M(M) & Left \& right & Ritz generalized singular vectors &  Small \& large
		\\\hline
	\end{tabular}
\end{table}

\subsection{MNIST Inpainting}\label{subsec:mnist}
\begin{figure}[!ht]
	\centering
	\includegraphics[width=\imgsize\linewidth]{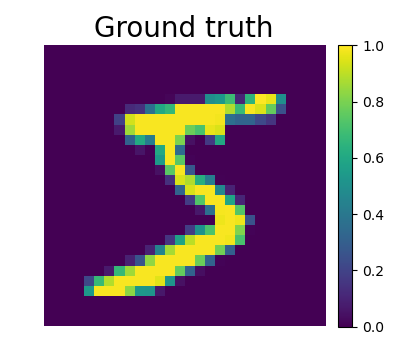}
	\includegraphics[width=\imgsize\linewidth]{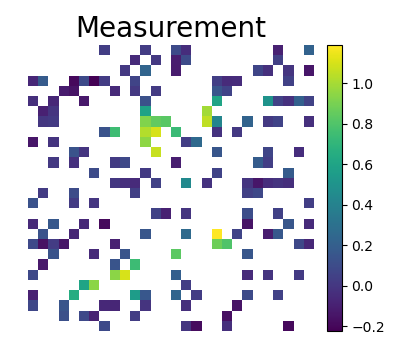}
	\includegraphics[width=\imgsize\linewidth]{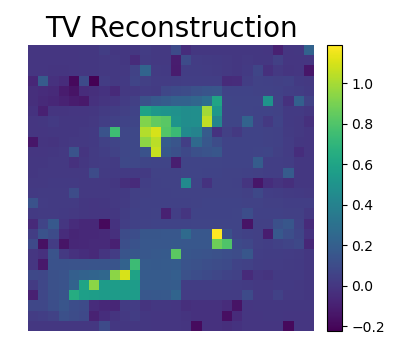}
	\includegraphics[width=\imgsize\linewidth]{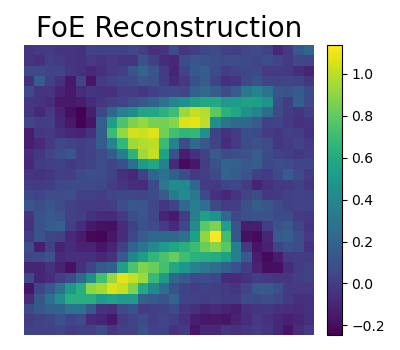}
	\caption{Data and reconstructions for the inpainting problem. We see that the Fields of Experts regularizer with optimized parameters can provide a much better reconstruction than TV.}\label{fig:inpaint}
\end{figure}

In order to numerically explore the effectiveness of the proposed method and approximations,  we consider a small bilevel problem consisting of a single data pair $(\xtrue,\ynoise)$  
where $\xtrue$ is a $28\times 28$ MNIST \cite{yannlecun2010mhd} image (so that $n=784$) and $\ynoise$ is obtained by applying a random subsampling matrix $A$ to $\xtrue$ (so that only $30\%$ of the original pixels are retained) and by adding Gaussian white noise of level ${\norm{\op \xtrue - \ynoise}_2}/{\norm{\op \xtrue}_2}=0.3.$
The ground truth and associated measurement are displayed in Figure~\ref{fig:inpaint}. 
We solve \eqref{eq:ul2ex} to learn $N=3$ filters of size $5\times 5$ (so that $p=78$) with $\phi(s)=s^2.$

We employ gradient descent with backtracking linesearch, as outlined in Algorithm~\ref{alg:hoag}, to solve \eqref{eq:ul2ex}. 
The lower level is solved via L-BFGS (history size $10$) with backtracking linesearch and we terminate when the gradient norm is smaller than $10^{-3}.$
We employ (R)MINRES to solve each Hessian system with a maximum of $500$ iterations where, unless specified otherwise, we terminate based on the residual norm $ \|{H^{(i)} w_k^{(i)} - g^{(i)}\|}_2 < \delta$  with a stopping criterion tolerance of $\delta = 10^{-2}.$ 
We find that early stopping is always achieved by (R)MINRES.
We consider a recycle dimension of ~$s=30,$ with this choice supported by one of the experiments which explores how the dimension of the recycle space affects performance.
Various stopping tolerances for (R)MINRES and L-BFGS were considered and the ones selected were such that the number of iterations of MINRES and L-BFGS were small while the determined solution of the bilevel problem achieved a hypergradient norm comparable to a high-accuracy run.

\begin{figure}[!ht]
	\centering
	\includegraphics[width=.8\linewidth]{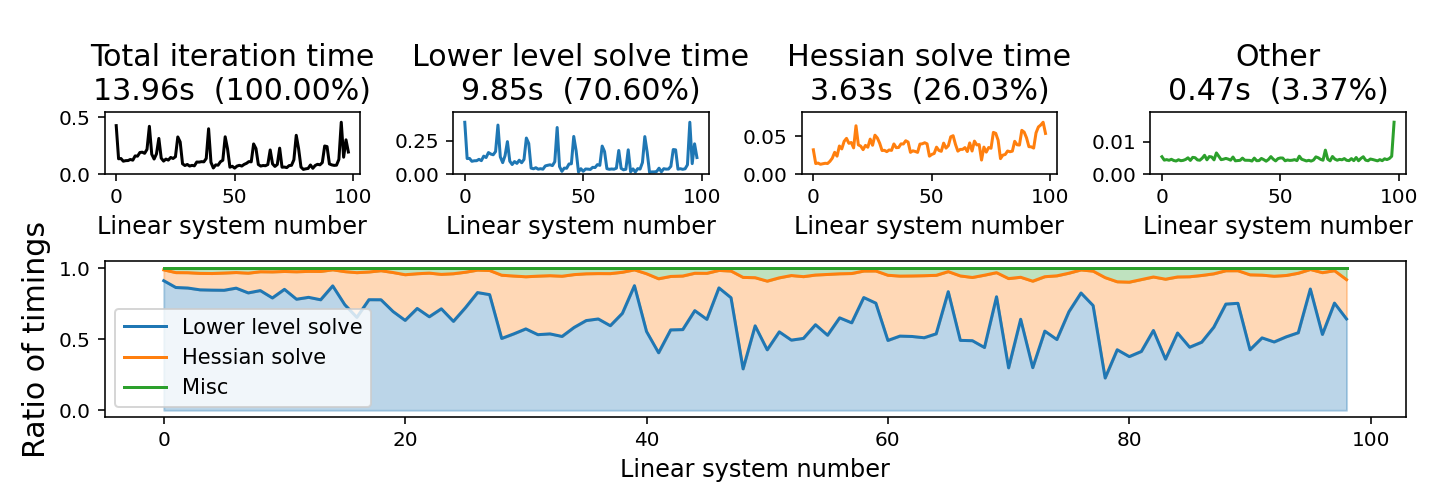}
	\caption{Breakdown of timings for solving the inpainting bilevel problem. 
		The titles of the subfigures in the top row indicate the total time (seconds) and percentage contribution towards the  total computation time.
		Each entry in the plot of timings for lower level solves is the sum of both the time for the initial  solve and timings of all solves performed within backtracking linesearch for that iteration of gradient descent for the upper level problem.
		Over 26\% of computation time is spent solving the sequence of Hessian systems.
	}\label{fig:og-timings}
\end{figure}

In Figure~\ref{fig:inpaint} we display reconstructions of the image used in~\eqref{eq:ul2ex}  obtained by Total Variation (TV) regularization, a common choice of regularizer for imaging applications, as well the FoE regularizer.
For the TV reconstruction, we selected the  regularization parameter achieving the smallest upper level cost among 20 logarithmically equispaced values between $10^{-16}$ and $10^1$. 
We see that due to the high undersampling, TV gives a poor reconstruction whereas, in learning filters specific to recovering that image, a learned regularizer provided by FoE recovers a meaningful reconstruction.
A breakdown of the computational timings of the bilevel solve is displayed in Figure~\ref{fig:og-timings}, where we recall that we employ backtracking linesearch to determine the step length for the upper level problem, and so solve more lower level problems than number of linear systems; see also Algorithms \ref{alg:backtrack} and \ref{alg:hoag}.  
We see that over 26\% of computation time is spent solving the sequence of Hessian systems.
Thus,  Krylov subspace recycling methods may provide a viable approach to reduce the computational cost by making such task 
more efficient.

\subsubsection*{Justification for recycling}
\begin{figure}
\centering
\includegraphics[width=\imgsize\linewidth]{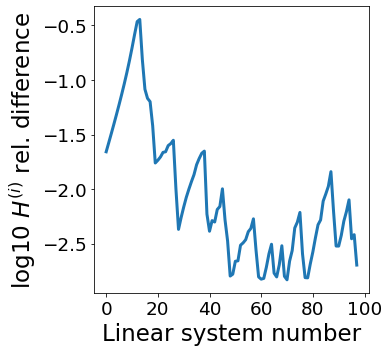}
\includegraphics[width=\imgsize\linewidth]{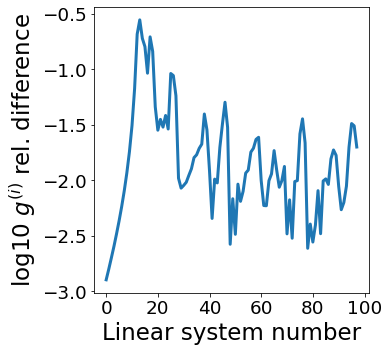}
\includegraphics[width=\imgsize\linewidth]{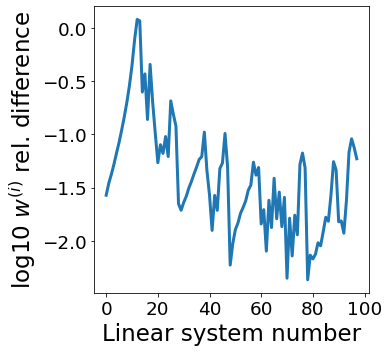}
\caption{Relative differences in the Frobenius norm for the sequence of, from left to right: the Hessians, right hand sides, and determined  solutions of the linear systems.
	Since relative differences of order $10^{-2}$ are observed, recycling may provide a viable approach to speed up computations.}\label{fig:justify}
	\end{figure}

We first experimentally verify that the linear systems encountered in the numerical solution of the inpainting bilevel problem are indeed similar. 
We consider the sequence of the Hessians $H^{(i)}$, right hand sides $g^{(i)}$ and solutions $w^{(i)}$ that were determined without recycling, and compute the relative difference in the Frobenius norm for each of these 3 sequences.
In Figure~\ref{fig:justify} we plot such quantities 
and we eventually see 
relative differences of order $10^{-2}$ 
in all instances.
Thus, Krylov subspace recycling may be a viable approach to reduce the overall computational cost. We emphasise that, in general, although outside the scope of this paper, such experimental validation may be supported by theoretical derivations, exploiting known properties of the Hessians $H^{(i)}$ (linked to the properties of the matrix $A$ and the filters $\widetilde{\theta}^{(j)}$), and how these are affected by multiplications by $J$, 
leading to an analysis similar to the one in \cite{kilmer2006rsi}.
Specifically, in \cite{kilmer2006rsi}, for a sequence of linear systems associated to the discretization of differential operators, and 
for consecutive systems that differ by a low enough rank plus small enough norm perturbation, the authors conclude that  approximate invariant subspaces (good for recycling)
are given by the (approximate) eigenvectors corresponding to the (approximate) small eigenvalues, that appear in disjoint enough clusters.
	
	\subsubsection*{MINRES vs CG}
	To employ recycling we require a method that explicitly  constructs a basis of the generated solution space and can be readily modified to incorporate the recycle space in the determined solution.
	For this reason we choose the recycling variant of MINRES, since the symmetry of each Hessian is still exploited to yield a computationally cheap method.
	This is in contrast to the standard choice of conjugate gradient (CG), where both the symmetry and positive definiteness of the Hessians are exploited.
	While MINRES has been empirically shown to be competitive to CG even for positive definite systems~\cite{chin2012cme}, we want to verify that, in using (R)MINRES, we are not just trying to recover the performance we would otherwise obtain if we employed CG.
	
\begin{figure}[!ht]
\centering
\includegraphics[width=\linewidth]{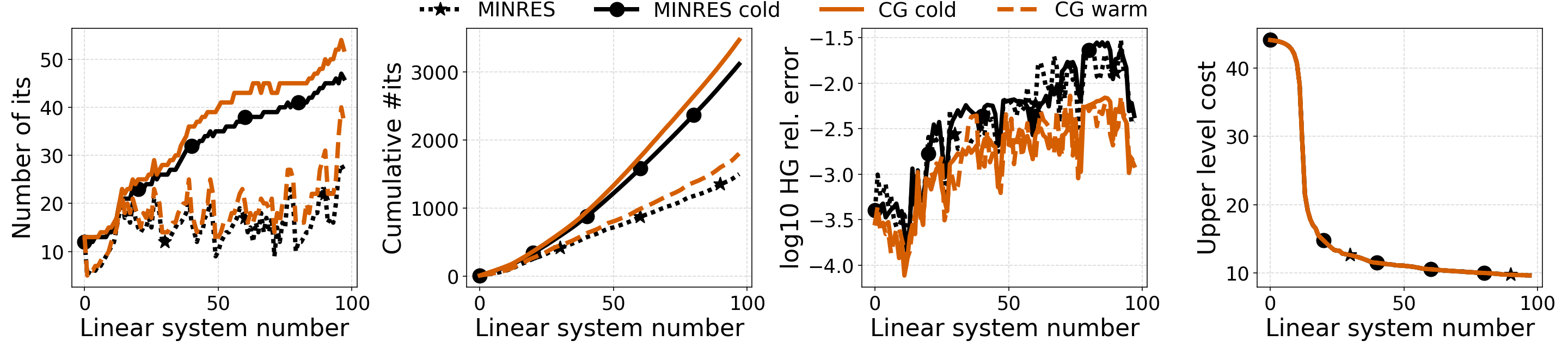}

\caption{Comparison between solving the sequence of Hessians without recycling using conjugate gradient (CG) and MINRES. Both methods have a 
	similar performance, with MINRES terminating earlier and the hypergradients associated to CG being slightly higher quality.
}\label{fig:cg-vs-minres}
\end{figure}

	In Figure~\ref{fig:cg-vs-minres} we plot, as a function of the index of the linear system sequence, the number of iterations, cumulative number of iterations, relative error of the associated hypergradients, and determined upper level cost for  MINRES and CG for both cold and warm starting.
	We see that MINRES and CG achieve similar performance for the considered problem.
	Employing warm starting yields a reduction in the number of iterations of both methods, but yields essentially identical quality of the determined  hypergradient.

\subsubsection*{Approximated versus full dimensional matrix decompositions}
We have motivated Ritz vectors and Ritz generalized singular vectors as approximations of eigenvectors and generalized singular vectors respectively, restricted to some subspace $W.$
Computing the full-dimensional matrix decompositions (i.e. taking $W=I$ in either \eqref{eq:ritz-matrix-decomp} or \eqref{eq:gsvd}) is not viable in practice due to the high computational cost. However, to provide insight into how the approximations compare to them, we solve the sequence of systems with Ritz-S, RGen-L(R), and the corresponding full-dimensional matrix decomposition strategies, Eig-S and GSVD-L(R), where we still select the recycle space to have dimension at most $s=30$.

\begin{figure}[!ht]
\centering
\includegraphics[width=\linewidth]{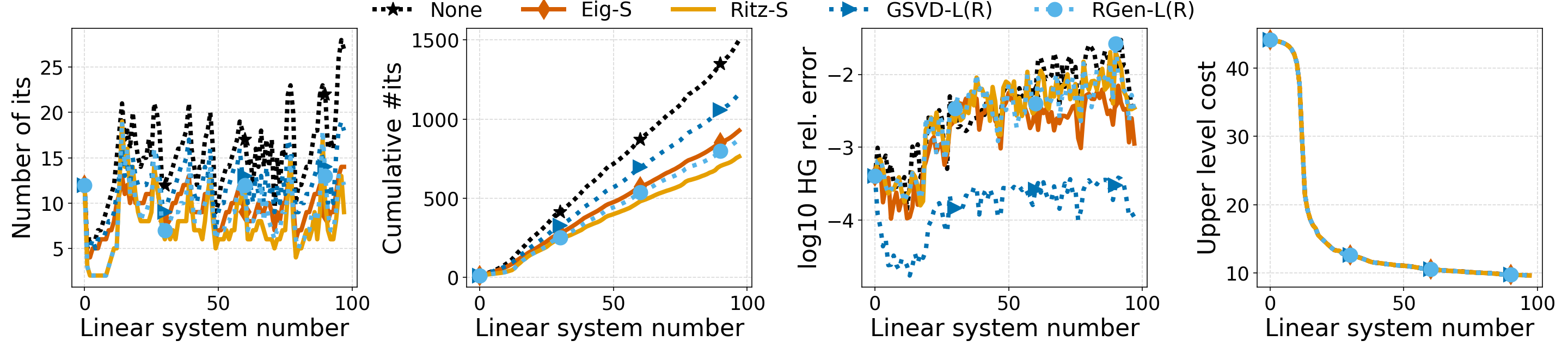}

\caption{Comparison between the performance of Ritz vectors to the eigenvectors they approximate, and Ritz generalized singular vectors to generalized singular vectors.
	A summary of the acronyms appearing in the legends is given in Table~\ref{tab:acro}.
}\label{fig:ritz-v-eig}
\end{figure}

In Figure~\ref{fig:ritz-v-eig} we plot the number of iterations to terminate RMINRES for each system, the cumulative number of RMINRES iterations, the relative error of the obtained hypergradients, and the associated upper level functional. 
We see that Ritz-S and RGen-L(R) regularly terminate earlier than their full-dimensional matrix decomposition counterparts, Eig-S and GSVD-L(R), respectively.
Since the residual norm stopping criterion is used in all of these strategies, this earlier termination of Ritz-S and RGen-L(R) suggests that considering $W\neq I$ can compact useful information of the solution space from the full-dimensional matrix decomposition. Indeed, Ritz vectors can be viewed as projected approximations of eigenvectors.
We see that GSVD-L(R) yields hypergradients an order of magnitude more accurate than the other methods, which obtain similar quality hypergradients to one-another. 
This supports the claim that encoding information of both the Hessian and adjoint Jacobian in the recycle space is a valid strategy to yield high quality hypergradients.

\subsubsection*{Small versus large Ritz 
values}
It is not always clear whether one should select vectors to recycle that correspond to the smallest or largest (harmonic) Ritz values, or possibly a mixture of the two; the same applies to Ritz generalized singular vectors. 
To explore this aspect, we solve the sequence of Hessian systems using all combinations of \textit{Size} and \textit{Side} for the Ritz, HRitz, and RGen recycling vector types.
For each \textit{Side} of the new recycling strategy, the \textit{Size} that yields the biggest reduction in iterations is employed again but using the new stopping criterion (NSC) \eqref{eq:hg-err-stop-lowrank} that approximates the hypergradient error.
We denote this by appending \textit{-NSC} to the acronym of the recycling strategy.

\begin{figure}
\centering
\includegraphics[width=\linewidth]{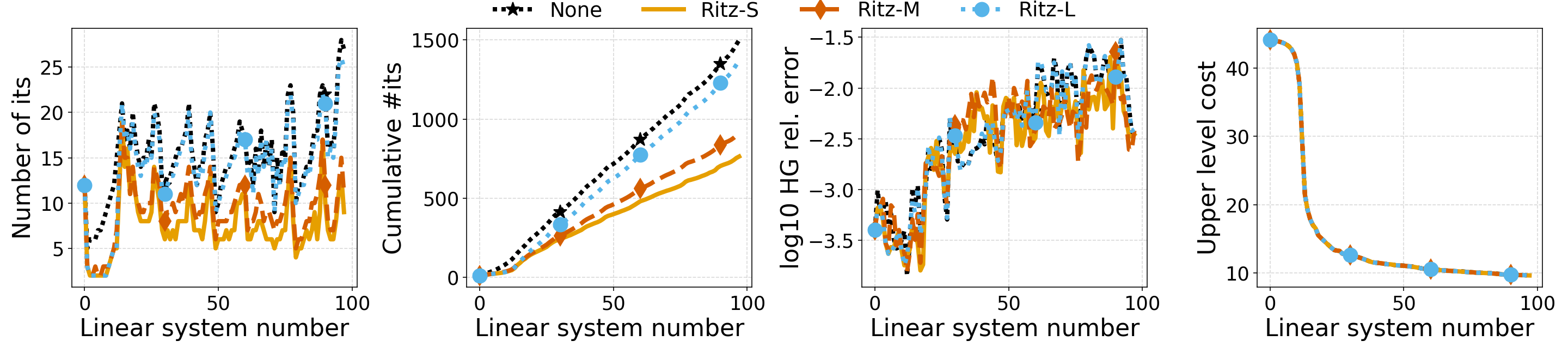}
\includegraphics[width=\linewidth]{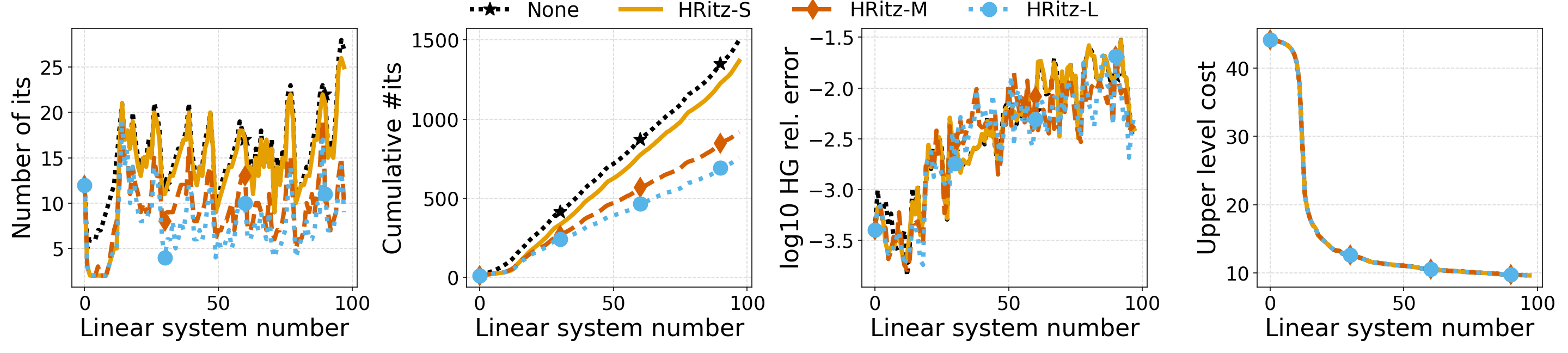}
\includegraphics[width=\linewidth]{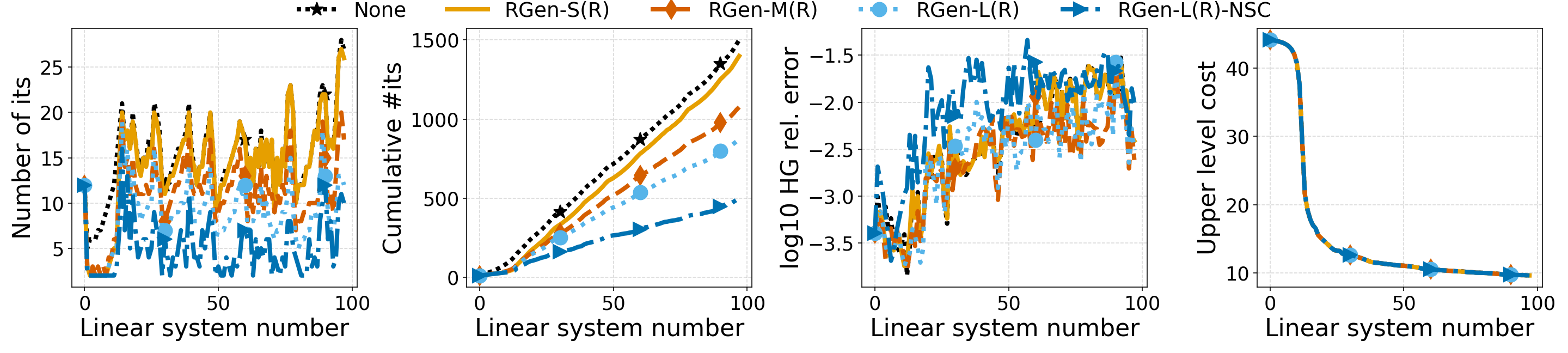}
\includegraphics[width=\linewidth]{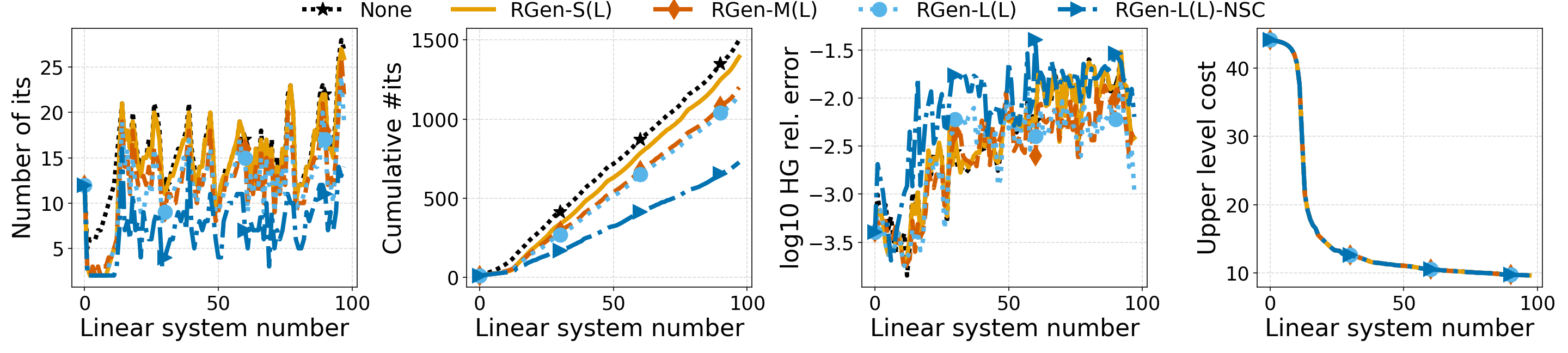}
\includegraphics[width=\linewidth]{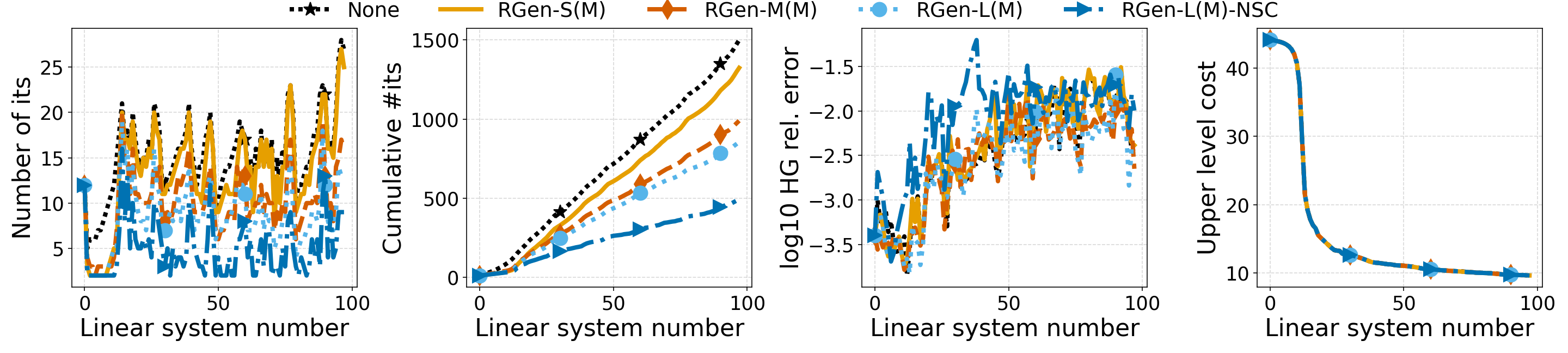}
\includegraphics[width=\linewidth]{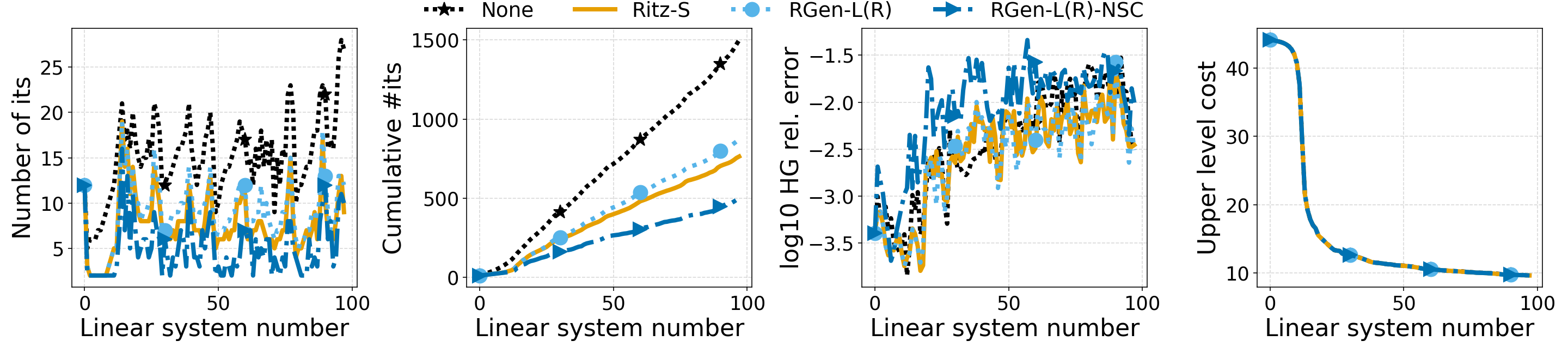}

\caption{Performance of solving the sequence of Hessian systems for the different recycling strategies. 
	From top to bottom row: Ritz vectors; right generalized singular Ritz vectors; harmonic Ritz vectors; left generalized singular Ritz vectors;
	mixture of left and right generalized singular Ritz vectors;    
	best performing recycle strategies that utilize Ritz vectors and Ritz generalized singular vectors.
	Employing the stopping criterion that approximates the hypergradient error yields a large reduction in iterations with no impact on the quality of the determined parameters.
	A summary of the acronyms appearing in the legends, located above each row, is given in Table~\ref{tab:acro}.}\label{fig:size-comp}
	\end{figure}

We plot the performance of the considered recycling strategies in Figure~\ref{fig:size-comp}.
We see that Ritz vectors associated with the smallest Ritz values (i.e., Ritz-S) yield the largest reduction in iterations, with those associated with large Ritz values achieving little improvement over not employing recycling at all. 
For harmonic Ritz vectors, information associated with the largest harmonic Ritz values results in the largest decrease of iterations, comparable to Ritz-S.
Similarly, we see that Ritz generalized singular vectors associated with the largest Ritz generalized singular values achieve the largest reduction in iterations. 
Given the performance of Ritz vectors this is to be expected. Indeed, should $\jacadj W=I,$ the Ritz generalized singular values are the reciprocals of the Ritz values.
Furthermore, using right rather than left Ritz generalized singular vectors yields a larger reduction in number of iterations. 
Utilizing the new stopping criterion \eqref{eq:hg-err-stop-lowrank} roughly halves the total number of iterations compared to using the residual norm.
Indeed, in the bottom row of Figure~\ref{fig:size-comp} we see that although Ritz vectors perform better than Ritz generalized vectors when using the residual norm, in employing the new stopping criterion, right Ritz generalized vectors yield a greater reduction in number of iterations.
We see that the new stopping criterion often results in a larger relative error in the associated hypergradient compared to other recycling strategies and stopping criteria. 
It should be noted that the methods which do not employ the new stopping criterion do not have a way to evaluate (an approximation of) the hypergradient error, which 
is unknown in practice. 
In contrast, the strategies that employ the new stopping criterion terminate when an approximation of the hypergradient error is $10^{-2},$ which, in this example, roughly coincides with the magnitude of the hypergradient relative error.

\subsubsection*{Iterations to reach a given hypergradient accuracy}
We have seen in the previous numerical experiment that, for the stopping criterion threshold of $\delta = 10^{-2},$ determined solutions of the sequence of Hessian systems yield hypergradients that achieve a relative error of around $10^{-2}.$ 
To determine how well recycling captures information relevant for the hypergradient, we use the true hypergradient  error \eqref{eq:hg-error-stop}  as the stopping criterion and solve the sequence of linear systems for Ritz-S, RGen-L(R), and no recycling.

\begin{figure}[!ht]
\centering
\includegraphics[width=\linewidth]{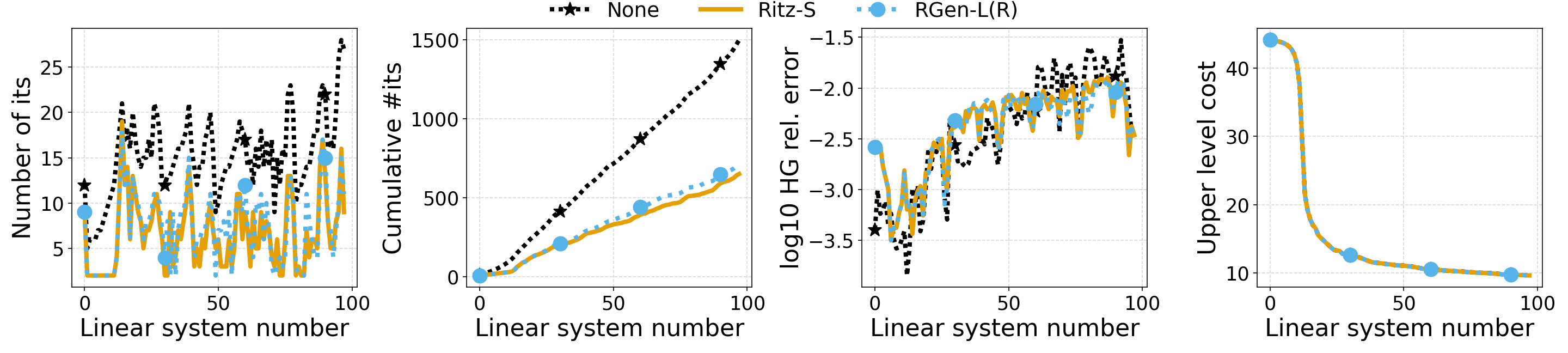}

\caption{Performance of solving the sequence of Hessian systems for different recycling strategies when the true hypergradient error is employed as the stopping criterion. 
	A summary of the acronyms appearing in the legends is given in Table~\ref{tab:acro}.}\label{fig:hg_as_stop_crit}
	\end{figure}

We display the results using the true hypergradient error  as the stopping criterion inFigure~\ref{fig:hg_as_stop_crit}.
We see that no recycling takes 1447 total iterations, which is little reduction compared to using a residual norm as a stopping criterion which required 1500 iterations (see Figure~\ref{fig:size-comp}).
This indicates that, in building the Krylov subspace from scratch, all iterations of MINRES are required to generate a solution space that can lead to a good hypergradient.
In contrast, by recycling information from the previous solution space,  Ritz-S and RGen-L(R) reduce from 764 and 871 total when using the residual norm (see Figure~\ref{fig:size-comp}) to 652 (81.9\%) and  713 (85.3\%), respectively, when using the true hyptergradient error.
This is to be expected since an accurate solution to the linear system is only a proxy for an accurate hypergradient. 
Interestingly, in Figure~\ref{fig:size-comp} we see that RGen-L(R)-NSC, which employed the approximation \eqref{eq:hg-err-stop-lowrank} of the hypergradient error as the stopping criterion, took 500 iterations. This indicates that the proposed approximation \eqref{eq:hg-err-stop-lowrank} often acts as a lower bound of the true hypergradient error \eqref{eq:hg-error-stop}.
We explore this aspect more thoroughly in the next experiment.

	\subsubsection*{Approximation of the hypergradient error}
	To understand how good our proposed approximation of the hypergradient error  \eqref{eq:hg-err-stop-lowrank} is, with respect to the true hypergradient error \eqref{eq:hg-error-stop}, 
	we solve the sequence of systems with RGen-L(R)-NSC, i.e. we employ the approximated hypergradient error \eqref{eq:hg-err-stop-lowrank} as the stopping criterion,
	and study the iterates from each RMINRES call with respect to quantities relevant to the proposed approximation, of which we give details below.
	
	Recall that the approximated hypergradient error \eqref{eq:hg-err-stop-lowrank}   requires computing the GSVD involving a projection of $H^{(i)}$ onto a subspace consisting of information related to $H^{(i-1)},$ namely, $W^{(i)} = [U^{(i-1)}, \; V^{(i-1)}].$ 
	Thus, errors in the approximation can occur not only in the fact that $W^{(i)}\neq I$, but also by potential suboptimality of the subspace spanned by $W^{(i)}.$
	
	Since we solve the sequence of Hessians with RGen-L(R)-NSC, we can evaluate the approximated hypergradient error that utilizes the GSVD of $H^{(i)}$ projected onto $W^{(i+1)} = [U^{(i)} ,\; V^{(i)}],$ that is, information associated with the current linear system, rather than the previous one. 
	Although this is only possible if the entire sequence is solved in advance, in this way, we can isolate the effect of both performing a projection, and performing a projection onto a subspace related to the previous system.

\begin{figure}[!ht]
\centering
\includegraphics[width=\imgsize\linewidth]{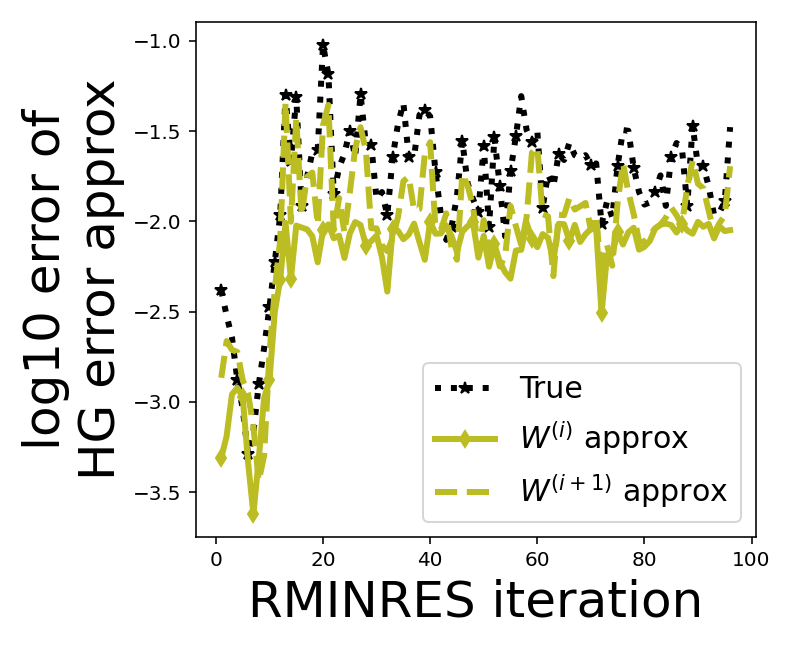}
\caption{Comparing the true hypergradient error to our approximation, with projection subspaces $W^{(i)}$ and $W^{(i+1)}$ associated to the previous and current systems, respectively, evaluated at the solution for each linear system solve. The sequence of Hessian systems is solved by RMINRES with RGen-L(R)-NSC.}\label{fig:approx_stop_all}
\end{figure}

	In Figure~\ref{fig:approx_stop_all} we plot evaluations of the true  hypergradient error, our proposed approximation that uses information of the previous system, and the approximation that utilizes information of the current system.
	The proposed approximation of the hypergradient error often serves as a lower bound for the other quantities, with the true hypergradient error often being an upper bound for the approximations. 

\begin{figure}[!ht]
\centering
\includegraphics[width=\imgsize\linewidth]{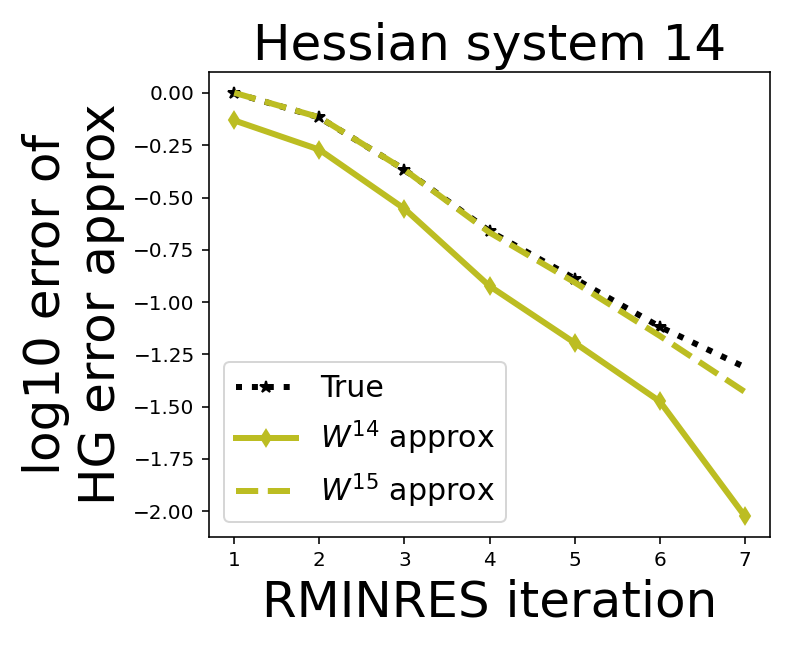}
\includegraphics[width=\imgsize\linewidth]{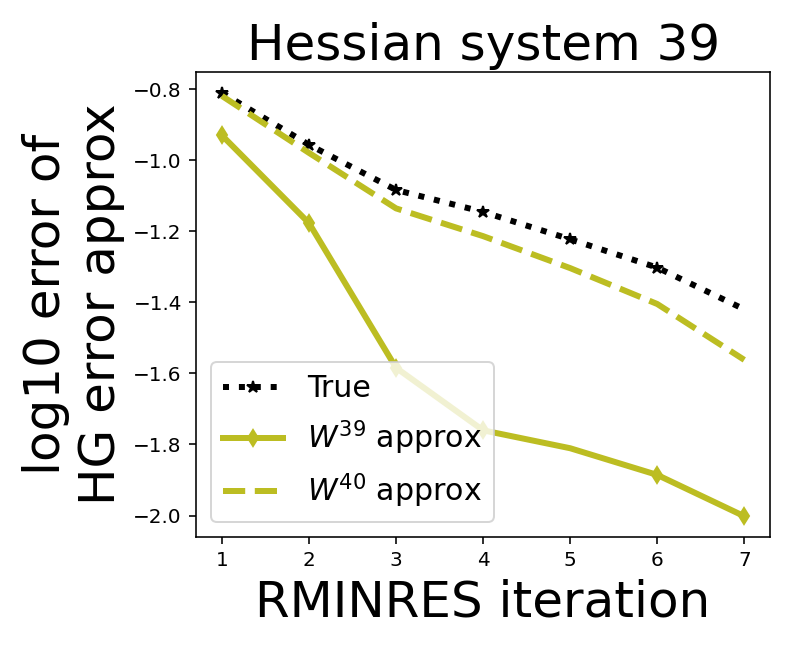}
\includegraphics[width=\imgsize\linewidth]{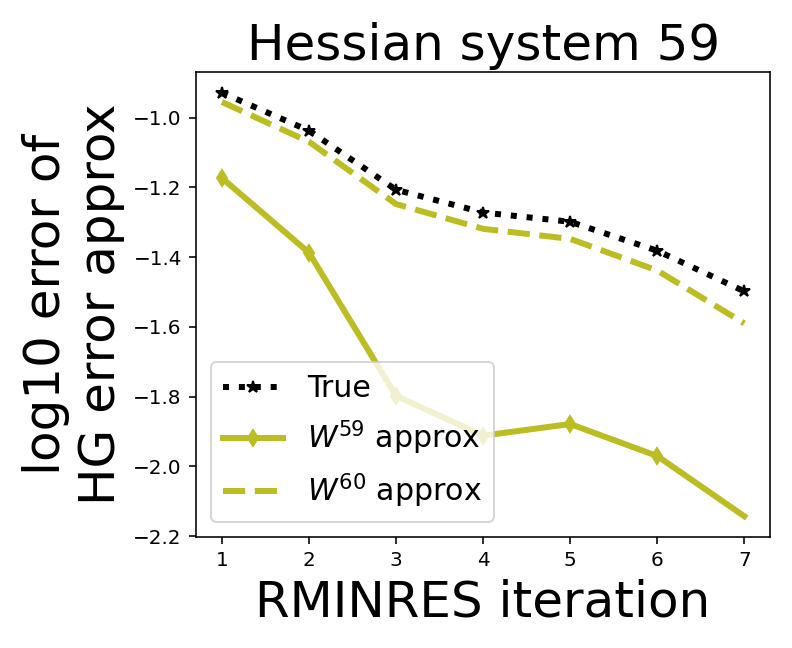}
\includegraphics[width=\imgsize\linewidth]{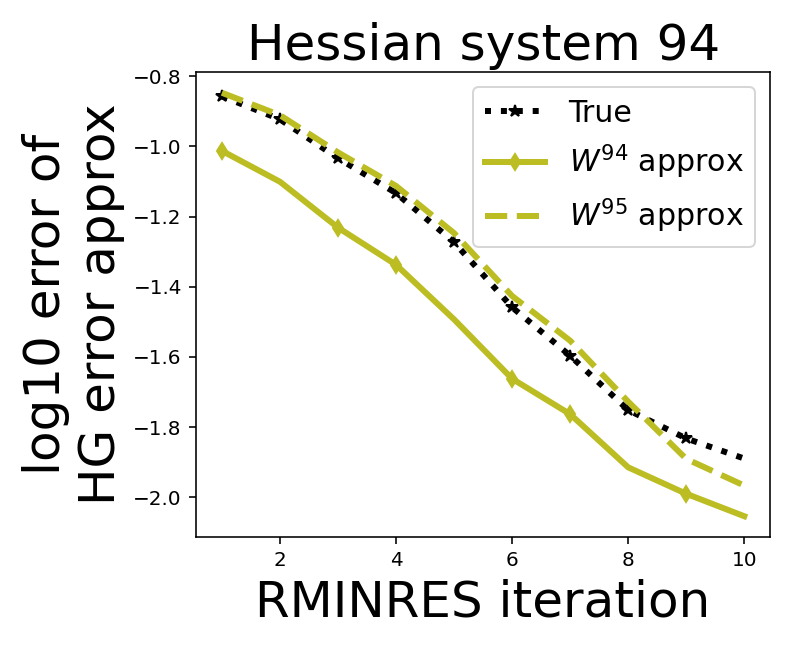}
\caption{Comparing the true hypergradient error to our approximation for different linear system solves. The sequence of Hessian systems is solved by RMINRES with RGen-L(R)-NSC. 
	Performing a projection onto the subspace $W^{(i+1)}$ associated with the linear system solution space provides a good approximation of the true hypergradient error, which indicates that the idea of performing a projection can provide a valid approximation. 
	The projection onto the solution space $W^{(i)}$ of the previous system, which is done in practice, provides a reasonable approximation, often serving as a lower bound.
}\label{fig:approx_stop}
\end{figure}

Next, we select a few different Hessian systems and consider how the different quantities change throughout each RMINRES run, and plot the results in Figure~\ref{fig:approx_stop}. 
We see that projection onto a subspace with information associated with the current linear system serves as a good approximation of the true hypergradient error, especially for early iterations. 
This could motivate using a strategy inspired by \cite{wang2007lto}, in which a recycle space to be employed for the next linear system solve is updated during the RMINRES run with information of the current system's solution space. 
In particular, one could consider a stopping criterion that employs this iteratively updated recycle space.
Similar to Figure~\ref{fig:approx_stop_all}, the proposed hypergradient error approximation  regularly serves as a lower bound for the other quantities.
Interestingly, for the 94th Hessian system one approximation is slightly larger than the true hypergradient error, which suggests that the proposed approximation may not be able to be formally viewed as a projected version. 

\subsubsection*{Inner versus outer recycling strategies}
To construct the $i$th recycle space $U^{(i)}$ which will be used in solving the linear system involving $H^{(i)},$ we currently consider a matrix decomposition (be it the eigendecomposition or GSVD) which involves $(W^{(i)})^T H^{(i)} W^{(i)},$
where $W^{(i)} = [U^{(i-1)},\; V^{(i-1)}]$ requires that we store the entire Krylov subspace basis between linear system solves.
In \cite{wang2007lto} it is proposed to consider a matrix decomposition concerning $(W^{(i)})^T H^{(i-1)} W^{(i)},$
as then one can avoid the need to store the whole basis and exploit matrix decompositions involving $H^{(i-1)}$ and $V^{(i-1)}.$ 
Specifically, for each linear system, our implementation requires storage of an $n\times(k_i+s_i)$ matrix, whereas \cite{wang2007lto} requires storage of a $n\times (\tilde k + s_i)$ matrix, where typically $\tilde k$ is chosen in advance such that $\tilde k \ll k_i$ for all $i;$ see  \cite{wang2007lto} for details.
We investigate the influence these two choices have on the recycling strategy.
We refer to a recycling strategy that employs the matrix proposed in \cite{wang2007lto}, that is, utilizing $H^{(i-1)},$ as an inner recycling strategy, since the recycle space for the next linear system is determined inside the RMINRES run for the current linear system.
Similarly, we refer to the choice that we primarily employ, utilizing $H^{(i)},$ as an outer recycling strategy.

\begin{figure}[!ht]
\centering
\includegraphics[width=\linewidth]{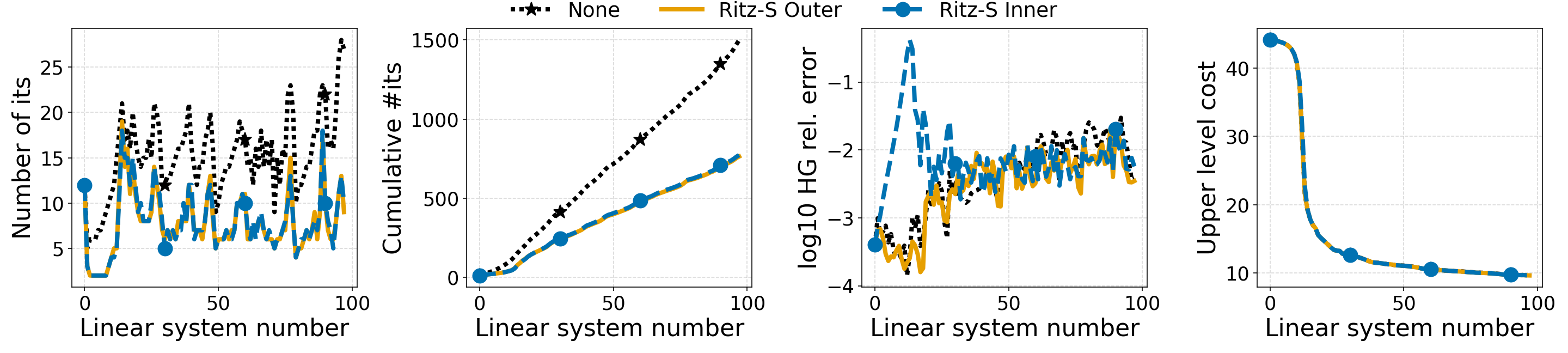}
\includegraphics[width=\linewidth]{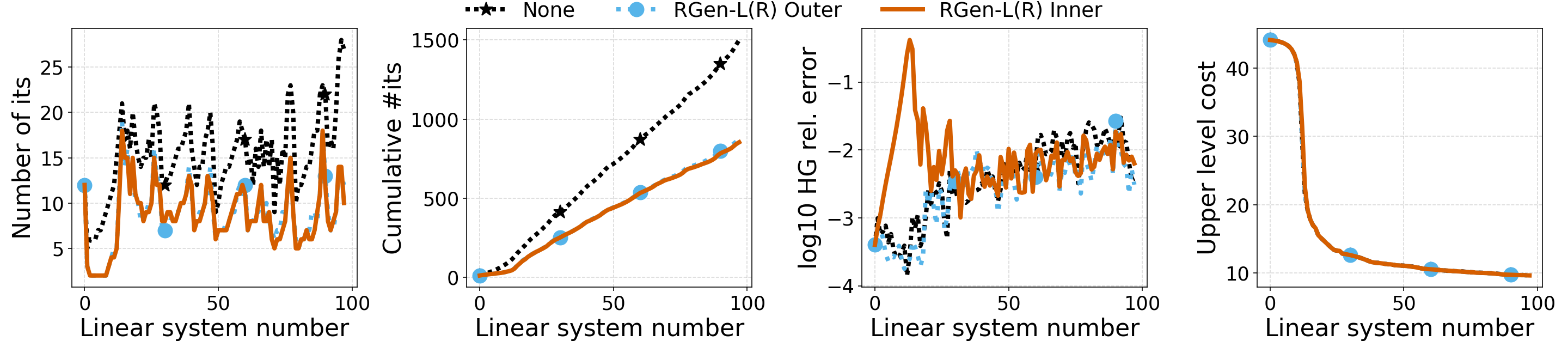}
\includegraphics[width=\linewidth]{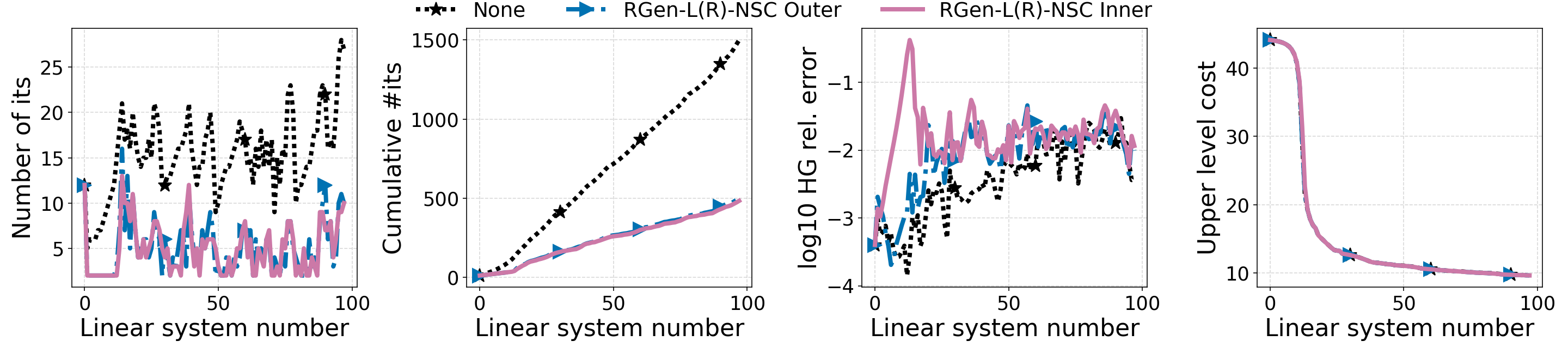}

\caption{Comparing the inner and outer variant of recycling using, from top to bottom row: 
	Ritz vectors associated with the smallest Ritz values; 
	left Ritz generalized singular vectors associated with the largest Ritz generalized singular values and the residual stopping criterion;
	left Ritz generalized singular vectors associated with the largest Ritz generalized singular values and new stopping criterion approximating the hypergradient error.
	The reduction in iterations is identical, but inner recycling yields worse quality hypergradients for the systems early in the sequence.
}\label{fig:inner-outer}
\end{figure}

In Figures~\ref{fig:inner-outer} we compare outer and inner variants of Ritz-S, RGen-L(R),and RGen-L(R)-NSC.
Both inner and outer methods yield the same performance in terms of the reduction in number of iterations.
Interestingly, inner recycling provide lower quality hypergradients for systems early in the sequence, but soon provide results comparable to outer recycling methods.
Despite the difference in hypergradient quality, the proposed parameter updates utilizing said hypergradients obtain comparable performance in terms of the upper level cost for both inner and outer recycling methods.

\subsubsection*{Dimension of the recycle space}
We now investigate how the dimension of the recycling space affects performance.
Indeed, to both construct Krylov vectors that are orthogonal to the recycle space, as in the augmented Lanczos relation \eqref{eq:arnoldi_aug_matrix}, and update the part of the solution that lies within the recycle space, extra FLOPS are required for a single iteration of RMINRES versus standard MINRES; see supplementary material.
Additionally, the cost of the matrix decomposition employed to determine the Ritz-type vectors will depend on the dimension of the recycle space.
Thus, one must balance the reduction of iterations to the computational overhead associated with a given recycling space dimension.

\begin{figure}[!ht]
\centering
\includegraphics[width=\linewidth]{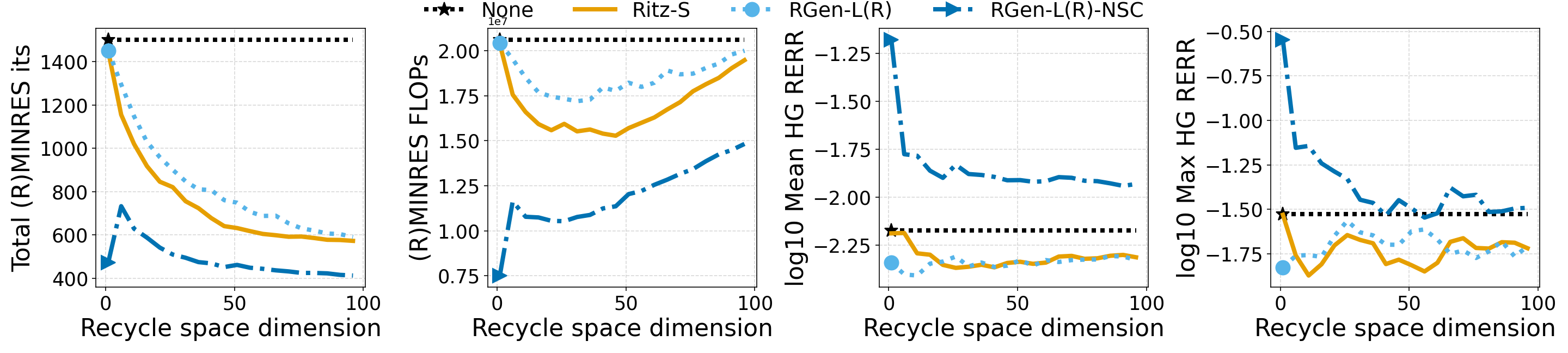}
\caption{Performance of Ritz-S and RGen-L(R) for different recycle dimensions.
	A summary of the acronyms appearing in the legends is given in Table~\ref{tab:acro}.
}\label{fig:dim-explore}
\end{figure}

For a given recycle space dimension and recycling strategy we can compute a sequence of hypergradient relative errors. 
In Figure~\ref{fig:dim-explore} we plot, as a function of the recycle space dimension, the maximum and mean of these hypergradient relative errors, as well as total number of RMINRES iterations, and total FLOPS of RMINRES.
Increasing the dimension reduces the total number of RMINRES iterations, since more information is passed between each linear system.
However, reduction in the number of iterations begins to stagnate for dimensions larger than 60. 
This is likely because the recycle space dimension is much larger than the number of iterations performed by RMINRES for a given system solve (around 10), and so entire solution spaces of several systems are fully acknowledged within the recycle space.
Thus, the recycle space becomes oversaturated, in that the inclusion of another Krylov solution space yields relatively little new information.

There is  a trade-off between the reduction of RMINRES iterations for a given recycling space dimension and the total number of FLOPS associated with all of the RMINRES runs. For the considered application,  a recycle space with dimension (at most) 30 provides a good tradeoff in this regard  - hence why we chose 30 to be the dimension of the recycle space throughout the  numerical experiments.
We see that the recycle space requires a dimension of at least 3 in order for the new stopping criterion to provide a meaningful approximation of the hypergradient error.
Indeed, RGen-L(R)-NSC terminates too early for a good quality hypergradient to be recovered for small recycle space dimensions.

\subsection{BSDS300 Deconvolution}\label{subsec:bsds}
\begin{figure}[!ht]
\centering
\includegraphics[width=.8\linewidth]{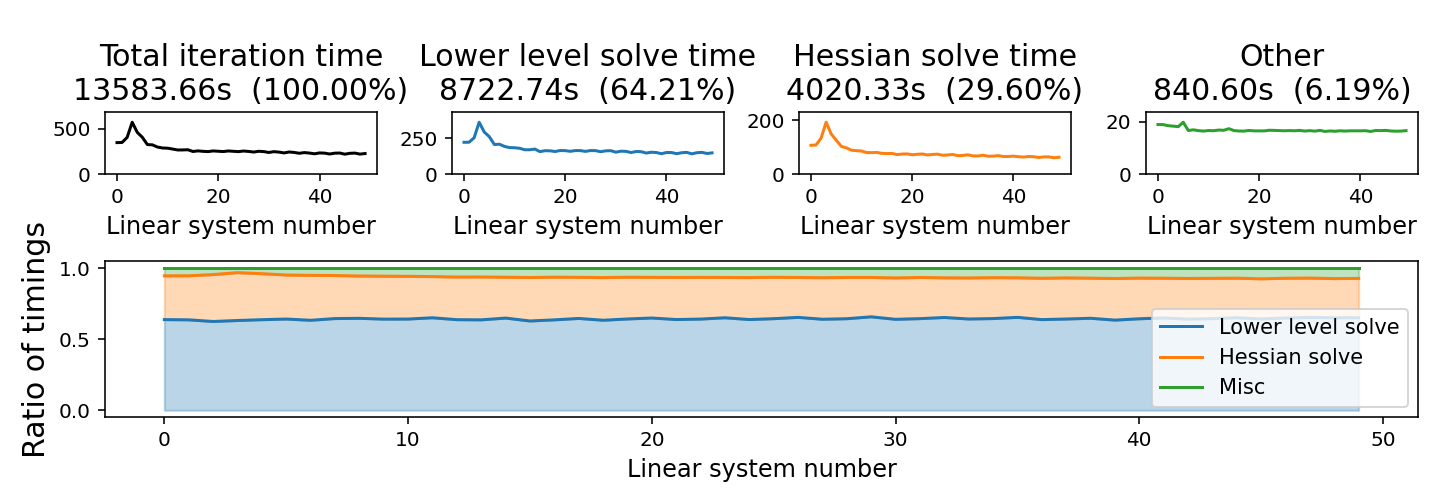}
\caption{Breakdown of timings for solving the deconvolution bilevel problem. 
	The titles of the subfigures in the top row indicate the total time (seconds) and percentage contribution towards the  total computation time.
	Over 29\% of computation time is spent solving the sequence of Hessian systems.}
\label{fig:bsds-timings}
\end{figure}

We now consider an application that resembles a  bilevel problem \eqref{eq:ul2ex} that is solved in practice.
For $64$ images  of the BSDS300 dataset \cite{MartinFTM01} we  take $8$ random $64\times 64$ crops (so $n=4096$ and $K=512$), which have been convolved with a Gaussian kernel of standard deviation $\sigma=3.0$ and then corrupted by Gaussian noise of level $0.2.$
We learn $N=24$ filters of size $5\times 5$ for the FoE regularizer~\eqref{eq:ul2ex} ($p=624$) with $\phi(s)=\log(1+s^2),$ which is motivated in the original work of \cite{roth2009fe}. 
Due to the non-convexity of $\phi,$ the Hessian of the lower level is no longer guaranteed to be positive definite - further justifying the choices of MINRES over CG for this application.
While the expression for hypergradients \eqref{eq:hg-full} may also no longer be valid, due to potential non-invertability of the Hessian, we follow common practice and construct hypergradients following the construction of \eqref{eq:hypgrad} regardless.
We solve \eqref{eq:ul2ex} without recycling using Adam for 50 epochs with a batch size of 16 and fixed steplength of $10^{-2}.$ 
For the lower-level we employ L-BFGS with backtracking (memory size 10).
A stopping tolerance of $\delta=10^{-3}$ is used for all involved optimizers and we allow for a maximum of 16000 iterations for both the lower level and Hessian system solvers. 
In line with~\cite{2017pocktnrd,kobler2023learninggraduallynonconveximage}, the filters are initialized with the non-constant 2D discrete cosine transform basis filters. Although the training data is split across mini-batches, each individual sample will, due to the fixed steplength, only require a single lower level and Hessian system solve each epoch.
We plot the accumulated timings of the lower level and Hessian system solvers across all training samples in Figure~\ref{fig:bsds-timings}, where we see that solving the Hessian systems takes over 29\% of computation time.

\begin{figure}[!ht]
\centering
\includegraphics[width=.8\linewidth]{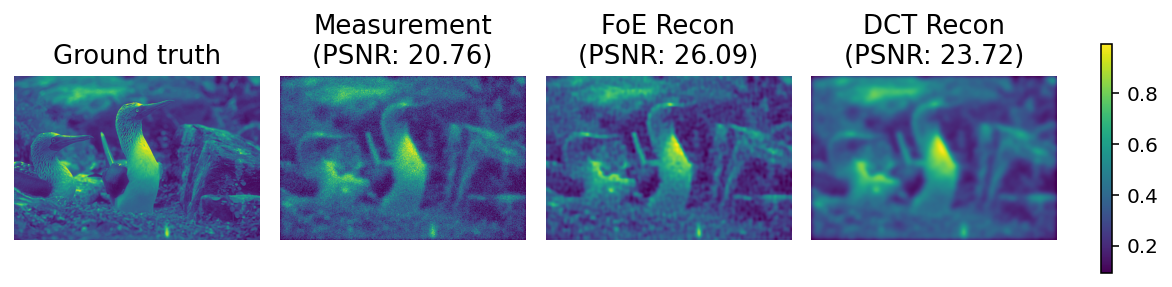}
\centering
\includegraphics[width=.35\linewidth]{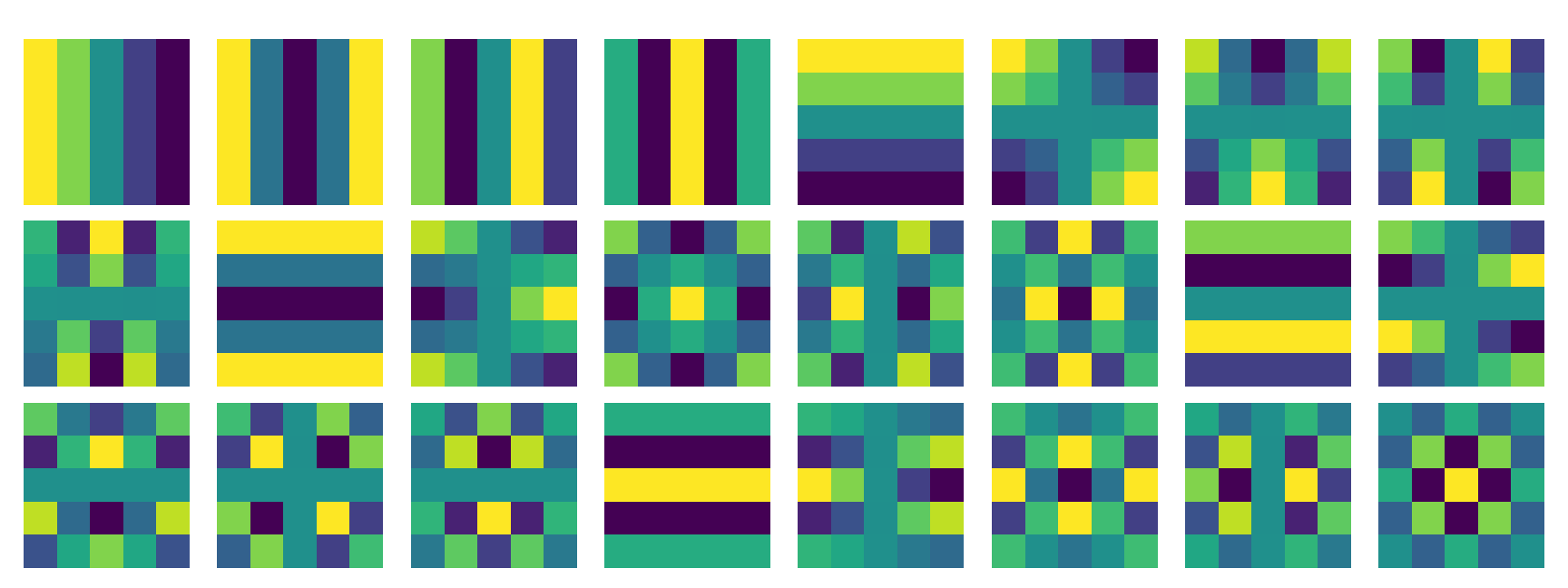}
\hspace{4ex}
\includegraphics[width=.35\linewidth]{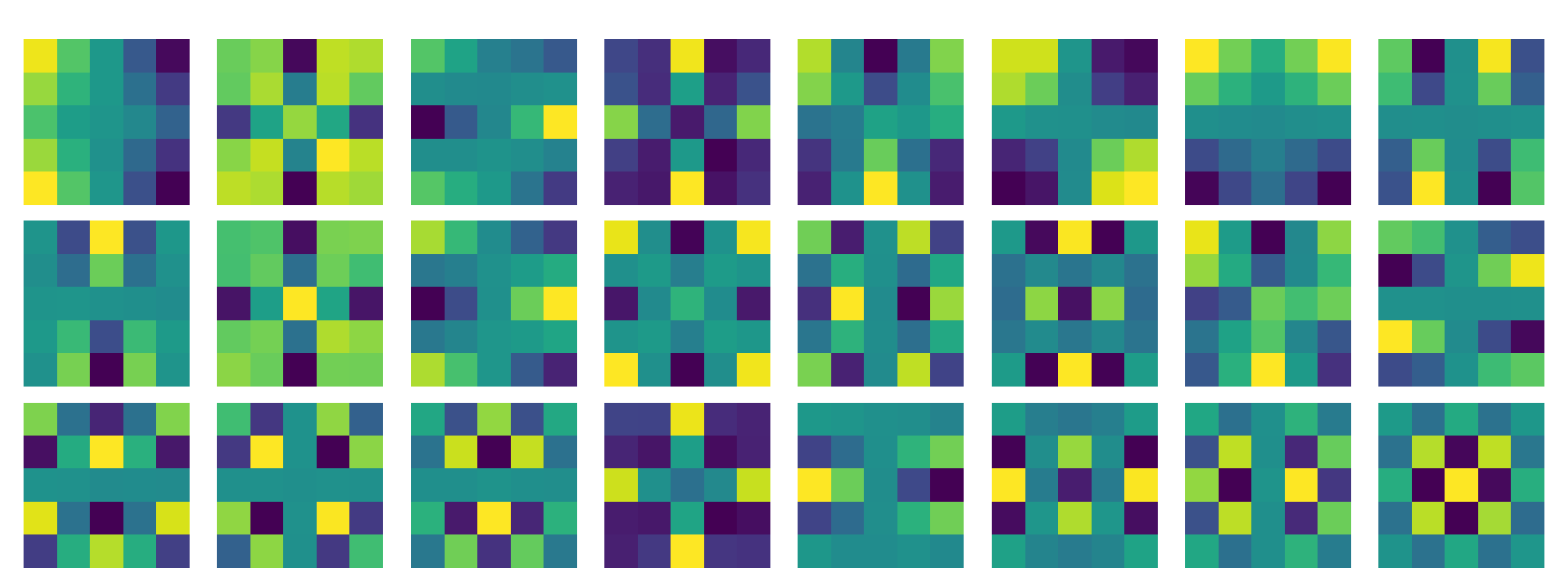}
\caption
{Top row: Data pair from the test dataset and associated reconstructions when using the learned FoE filters and the standard DCT basis. All pixel values are clipped to the interval $[0,1]$ when displayed.
	Bottom left: The initialized filters for the FoE model. 
	Bottom right: Filters after training. The filters associated with low frequencies have changed the most.
	The colours have been scaled to each individual filter.}
\label{fig:bsds-filters}
\end{figure}

Reconstructions of on an image from the test when using the learned and initial filters, that is, the DCT basis, are shown in top row of Figure~\ref{fig:bsds-filters}. 
The learned filters succesfuly recover a meaningful reconstruction. The initial and learned filters are displayed in the bottom row of Figure~\ref{fig:bsds-filters}, where we observe that the low frequency filters have changed the most.

\begin{figure}[!ht]
\centering
\includegraphics[width=.7\linewidth]{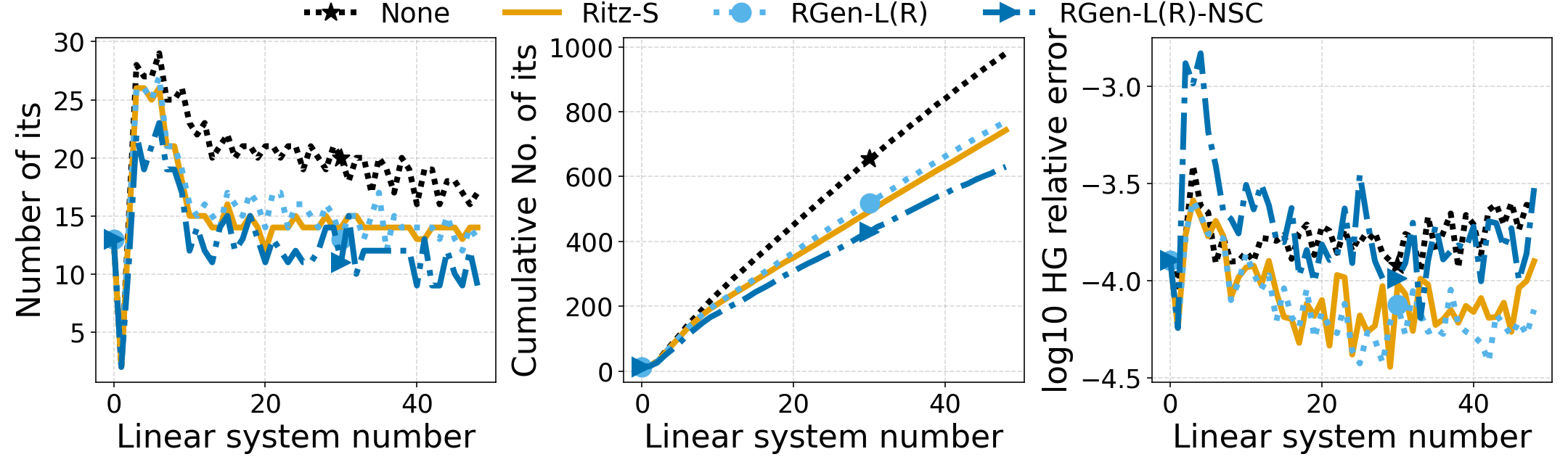}
\includegraphics[width=.7\linewidth]{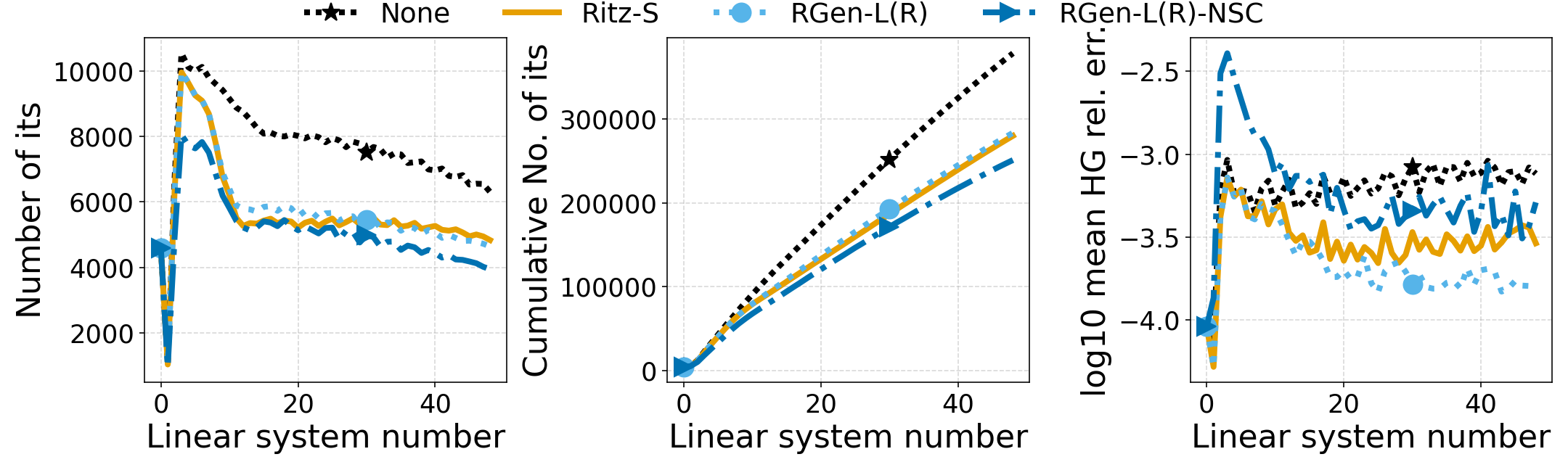}
\caption{Performance of solving the sequence of Hessian systems with the best performing strategies of Section~\ref{subsec:mnist}.
	Top row: results associated with the first training data sample. 
	Bottom row: accumulated performance for all training data.
	Employing the stopping criterion that approximates the hypergradient error yields the largest reduction in iterations with no impact on the quality of the determined parameters.
	A summary of the acronyms appearing in the legend is given in Table~\ref{tab:acro}. }
\label{fig:bsds-compare}
\end{figure}

We re-solve the sequence of Hessians associated with every sample of training data data using the recycling strategies that performed best in Section~\ref{subsec:mnist}, namely, {Rits-S}, {RGen-L(R)}, and {RGen-L(R)-NSC}. 
The performance of said strategies for solving the sequence of systems associated with first sample of training data is displayed in the top row of Figure~\ref{fig:bsds-compare}. 
Both Ritz-S and RGen-L(R) reduce the overall number of RMINRES iterations by approximately 27\%, with Ritz-S slightly out-performing RGen-L(R). 
Employing the new stopping criterion \eqref{eq:hg-err-stop-lowrank} leads to a reduction of 39\% while maintaining hypergradients of similar quality.
The bottom row of Figure~\ref{fig:bsds-compare} shows the accumulated performance across all sequences of systems solved during training. 
Both Ritz-S and RGen-L(R) yield a reduction of $26\%,$ while the RGen-L(R) obtains a larger reduction of $34\%.$ 
All methods attain a smaller reduction than what was observed in Section~\ref{subsec:mnist}. This is likely a result of using mini-batches in training, as potentially multiple parameter updates have been performed before a  specific training sample is encountered again.


\section{Conclusions and outlook}\label{sec:rec-conclusion}
To employ a gradient based method to solve a bilevel optimization problem, hypergradients must be computed by solving a sequence of lower level problems and Hessian systems, which leads to a high computational cost. 
To alleviate this cost, in this work we employed Krylov subspace recycling methods to solve the encountered sequence of Hessian systems,
which, due to the similarity between each system, results in a large reduction in the number of overall  iterations. 
We propose a novel recycling strategy based on the GSVD and a new concept, Ritz generalized singular vectors (associated to the GSVD of a couple of projected matrices), in order to acknowledge the bilevel setting in which the sequence of Hessians arise from.
This is in contrast to existing recycling strategies, where the recycle space is determined in isolation to how the solution of the linear system is to be utilized. 
Furthermore, by employing the quantities appearing in the projected GSVD, we can employ a stopping criterion for the Hessian system which approximates the hypergradient error - a quantity usually completely unavailable. 
We empirically show that the joint application of the new recycling strategy and stopping criterion leads to computations of hypergradients and hyperparameters of the same accuracy as the ones obtained using traditional solvers, with significant computational savings. 

Future work includes a rigorous analysis of both  optimality of Ritz generalized singular vectors and the accuracy of the hypergradient error approximation, as well as the study of alternative optimal approaches that take into account the premultiplication of $J$ in the hypergradient computations.
While extensive numerical experiments verify the proposed recycling strategy as a proof of concept, an implementation in a spirit similar to \cite{wang2007lto} is to be determined, wherein the Lanczos relation is exploited within the GSVD computation.
Large-scale applications remain to be explored, specifically, those in which solutions of linear systems are employed in a matrix-vector product, for example, differentiable architecture search \cite{ma2024mrd}. 
Using the proposed recycling strategy in conjunction with other computational saving techniques, such as inexact methods \cite{ehrhardt2021ido,ehrhardt2023aih,pedregosa2016hoa,salehi2024aif}, should be investigated, with special attention given to potential interaction regarding the hypergradient error approximation.

\bibliographystyle{siamplain}

\newpage
\appendix

\section{Calculating the residual vector for RMINRES}

To solve a linear system $H w_\star = g$ we employ  RMINRES \cite{wang2007lto}, detailed in Algorithm~\ref{alg:rminres}.
\begin{algorithm}
	\caption{RMINRES}\label{alg:rminres}
	\begin{algorithmic}[1]
		\STATE Initial guess $w_0\in\real^n,$ stopping residual tolerance $\epsilon>0,$ matrices $H\in\real^{n\times n},$ $U,C\in\real^{n\times s},$ and right hand side $g\in\real^n$
		\STATE $r_0 := g - Hw_0$
		\STATE $ w_0 :=  w_0 + UC^Tr_0$
		\STATE $r_0 := r_0 - CC^Tr_0$
		\STATE $\beta_1:= \norm{r_0};\quad v_1 := r_0 /\beta_1; \quad \zeta_0 := \beta_1$
		
		\STATE $(\tilde v_0,\tilde v_{-1} , \tilde b_0 , \tilde b_{-1}) :=  (0,0,0,0)\in\real^{n\times 4}$
		\STATE $(c_0,s_0,c_{-1},s_{-1}) := (1,0,0,0)\in\real^4$
		
		\FOR{k=1,...}
		
		\STATE Determine new Lanczos vector:
		
		\STATE $v_{k+1}:= H v_k;\quad b_k := UC^Tv_{k+1}$  \label{code:Hvk}
		\STATE $v_{k+1} := v_{k+1} - CC^T  v_{k+1}; \quad  \alpha_k := \langle v_k, v_{k+1}\rangle$ \label{code:vk_update}
		
		\STATE $v_{k+1} := v_{k+1} - \alpha_k v_k - \beta_k  v_{k-1}
		$
		\STATE $\beta_{k+1} := \norm{v_{k+1}};\quad v_{k+1}:= v_{k+1}/\beta_{k+1}$
		
		\STATE
		
		Apply previous two Givens rotations to new column of the tridiagonal matrix:
		\STATE$		\gamma_k := s_{k-2} \beta_k$
		\STATE $
		\delta_k := c_{k-1}c_{k-2} \beta_{k} + s_{k-1}\alpha_k$
		\STATE $\epsilon_k := -s_{k-1}c_{k-2} \beta_{k} + c_{k-1}\alpha_k$
		
		\STATE 		Specify and apply new Givens rotation to annihilate $\beta_{k+1}:$
		
		\STATE $\tilde\epsilon_k := \sqrt{\epsilon_k^2 + \beta_{k+1}^2}$
		\STATE $s_k := \beta_{k+1}/\tilde\epsilon_k;\quad c_k:= \epsilon_k/\tilde\epsilon_k$
		\STATE $\zeta_k := -s_k\zeta_{k-1}$
		
		\STATE Update solution:
		
		\STATE $\tilde v_k :=  (v_k - \gamma_k\tilde v_{k-2} - \delta_k \tilde v_{k-1})/\tilde\epsilon_k;\quad \tilde b_k :=  (b_k - \gamma_k\tilde b_{k-2} - \delta_k \tilde b_{k-1})/\tilde\epsilon_k$ \label{code:vb_update}
		
		\STATE $w_k := w_{k-1}  + c_k\zeta_{k-1} (\tilde v_{k} + \tilde b_k)$ \label{code:w_update}
		
		\STATE Stopping criterion
		
		\IF{$\abs{\zeta_k} < \epsilon$}
		\STATE \textbf{break}
		\ENDIF
		\ENDFOR
		\RETURN $w_k$
	\end{algorithmic}
	\caption{RMINRES, a recycling variant of MINRES, see \cite{wang2007lto} for details.}
\end{algorithm}
While MINRES provides an efficient way to determine the residual norm, to employ the alternative stopping criterion discussed in Section~\ref{sec:stopcrit} we require explicit access to the residual vector $r_k := g - Hw_k.$  
Towards this, notice that from Algorithm~\ref{alg:rminres} line~\ref{code:w_update} the residual vector  will satisfy
\begin{equation*}
	r_k = r_{k-1} - c_k\zeta_{k-1} H(\tilde v_k - \tilde b_k).
\end{equation*}
Thus we require efficient ways to compute both $H \tilde v_k$ and $H \tilde b_k.$
To compute $H\tilde v_k,$ it follows from line~\ref{code:vb_update} that
\begin{equation*}
	H\tilde v_k = \frac{1}{\tilde\epsilon_k} \left( Hv_k - \delta_k H\tilde v_{k-1} - \gamma_k \tilde  v_{k-2}  \right)
\end{equation*}
and we already compute $Hv_k$ in line~\ref{code:Hvk} for the construction of $v_{k+1}.$ 
Since $\tilde v_{k-2}$ is already stored in RMINRES, we need only store two additional vectors, $Hv_k$ and $H\tilde v_{k-1},$ to determine $H\tilde v_k.$ 

Similarly, to compute $H\tilde b_k,$ it follows from   Algorithm~\ref{alg:rminres} lines~\ref{code:vb_update} and \ref{code:Hvk} that   
\begin{align*}
	H\tilde b_k 
	&= \frac{1}{\tilde\epsilon_k} \left( Hb_k - \delta_k H\tilde b_{k-1} - \gamma_k \tilde  b_{k-2}  \right)
	\\
	&= \frac{1}{\tilde\epsilon_k} \left( HUC^THv_k - \delta_k H\tilde b_{k-1} - \gamma_k \tilde  b_{k-2}  \right)
	\\
	&= \frac{1}{\tilde\epsilon_k} \left( CC^THv_k - \delta_k H\tilde b_{k-1} - \gamma_k \tilde  b_{k-2}  \right),
\end{align*}
where the last equality follows from $C:=HU.$ 
Now $CC^THv_k$ is already computed in Algorithm~\ref{alg:rminres} line~\ref{code:vk_update} in order to ensure  $v_{k+1}$ is orthogonal to $C.$
Since $\tilde b_{k-2}$ is already stored in RMINRES, we need only store two additional vectors, $CC^THv_k$ and $H\tilde b_{k-1}$ to determine $H\tilde b_k.$

Thus, we can iteratively update the residual vector if we store 5 additional vectors in $\real^n$ (the residual vector itself and the 4 vectors discussed above).

\section{Computational cost of RMINRES}
While RMINRES can yield a reduction in the number of iterations compared to standard MINRES, each iteration of RMINRES requires more FLOPS; in order to update the part of the solution represented in the recycle space as well as to ensure orthogonality of each Lanczos vector to $C:=HU.$

We quickly recall that for $A\in\real^{n\times s}, B \in\real^{m\times n}, x\in\real^n,$ the number of FLOPS required to compute products $Ax$ and $AB$ is $n(2s-1)$ and $ms(2n-1)$ respectively.

A main cost will be the application of the Hessian $H\in\real^{n\times n},$ which has a naive matrix-vector product cost of $2n^2-n.$ 
However, depending on the choice of lower level cost function, it may be possible to determine efficient implementations of matrix-vector products with the Hessian.
For this reason, we will refer to the general cost of an application of $H$ as $H_{\text{cost}}.$

To determine the cost of MINRES, we may ignore products and terms that arise from applications of $U$ and $C$ in Algorithm~\ref{alg:rminres}.
Thus, without any recycling, we see that calculations before the main loop require  $H_{\text{cost}} + 4n$ FLOPS.
Within the for loop, construction of  $v_{k+1}$ costs  $H_{\text{cost}} + 7n$ and updating the solution costs 
$9n.$ 
Thus,  $k_{\text{total}}$ iterations of MINRES costs
\begin{equation}
	H_{\text{cost}} + 4n + k_{\text{total}}\left(H_{\text{cost}} + 16n \right)\label{eq:minres_cost}
\end{equation}
FLOPS.

Employing RMINRES involves additional computations, namely, the orthogonalisation of Lanczos vectors against $C,$ and updating the part of the solution that resides in the recycle space.
Notice that we can reuse the value of $C^Tr_0$ (and $C^Tv_{k+1}$) to avoid an unnecessary application of $C^T$ both in and outside the for loop.
With this in mind, the additional cost of RMINRES outside the  loop is $6ns.$ Thus, for every loop iteration, constructing $v_{k+1}$ and $b_k$ incurs an additional cost of $6ns-n-s$ and updating the solution an extra cost of $6n.$ 
It follows that  $k_{\text{total}}$ iterations of RMINRES costs
\begin{equation}
	H_{\text{cost}} + 4n + 6ns +  k_{\text{total}}\left(H_{\text{cost}} + 21n + 6ns  - s \right).\label{eq:rminres_cost}
\end{equation}
We remark that, due to the additional 3 term recurrence that is employed to update the solution with respect to the recycle space, each iteration of RMINRES has an additional cost independent of the recycle dimension. 
Finally,  \eqref{eq:rminres_cost} is only the cost associated with a run of RMINRES and neglects the cost to determine the recycle space itself. 
{We remark that, depending on the application (that is, the value of $H_{\text{cost}})$ and number of iterations RMINRES is expected to take, it is possible implement a more computationally effective algorithm than \ref{alg:rminres} \cite{MELLO20103101}. }

\end{document}